\documentclass[12pt]{amsart}
\usepackage{graphicx}
\usepackage{graphicx}
\usepackage{epsfig}
\usepackage{pstricks}

\newtheorem{thm}{Theorem}[section]

\newtheorem{defn}[thm]{Definition}
\newtheorem{lemma}[thm]{Lemma}
\newtheorem{conj}[thm]{Conjecture}
\newtheorem{cor}[thm]{Corollary}
\newtheorem{remark}[thm]{Remark}

\usepackage{amsmath}
\usepackage{amsxtra}
\usepackage{amscd}
\usepackage{amsthm}
\usepackage{amsfonts}
\usepackage{amssymb}
\usepackage{eucal}
\newcommand{\bmb}{\left( \begin{array}{rr}}
\newcommand{\enm}{\end{array}\right)}

\newcommand{\cQ}{\mathcal Q}
\newcommand{\cT}{\mathcal T}

\newcommand{\cP}{\mathcal P}

\newcommand{\C}{{\mathbb C}}
\newcommand{\Z}{{\mathbb Z}}

\newcommand{\bF}{{\mathbf F}}
\newcommand{\bG}{{\mathbf G}}

\newcommand{\bx}{{\mathbf x}}

\newcommand{\bt}{{\mathbf t}}
\newcommand{\bu}{{\mathbf u}}
\newcommand{\bv}{{\mathbf v}}
\newcommand{\bw}{{\mathbf w}}

\newcommand{\bR}{{\mathbf R}}

\newcommand{\bC}{{\mathbf C}}

\newcommand{\bz}{{\mathbf z}}

\newcommand{\al}{{\alpha}}

\numberwithin{equation}{section}

\begin{document}

\title{Twenty Vertex model and domino tilings of the Aztec triangle}
\author{Philippe Di Francesco} 
\address{Department of Mathematics, University of Illinois, Urbana, IL 61821, U.S.A. 
and \break
Institut de physique th\'eorique, Universit\'e Paris Saclay, 
CEA, CNRS, F-91191 Gif-sur-Yvette, FRANCE\hfill
\break  e-mail: philippe@illinois.edu
}

\begin{abstract}
We show that the number of configurations of the 20 Vertex model on certain domains with domain wall type boundary conditions is equal to the number of domino tilings of Aztec-like triangles, proving a conjecture from \cite{DG19}. The result is based on the integrability of the 20 Vertex model and uses a connection to the U-turn boundary 6 Vertex model to re-express the number of 20 Vertex configurations as a simple determinant, which is then related to a Lindstr\"om-Gessel-Viennot determinant for the domino tiling problem. The common number of configurations is conjectured to be
$2^{n(n-1)/2}\prod_{j=0}^{n-1}\frac{(4j+2)!}{(n+2j+1)!}=1, 4, 60, 3328, 678912...$ The enumeration result is extended to include refinements of both numbers.  
\end{abstract}

\maketitle
\date{\today}
\tableofcontents

\section{Introduction}

\subsection{Triangular Ice model combinatorics}

\begin{figure}
\begin{center}
\includegraphics[width=14cm]{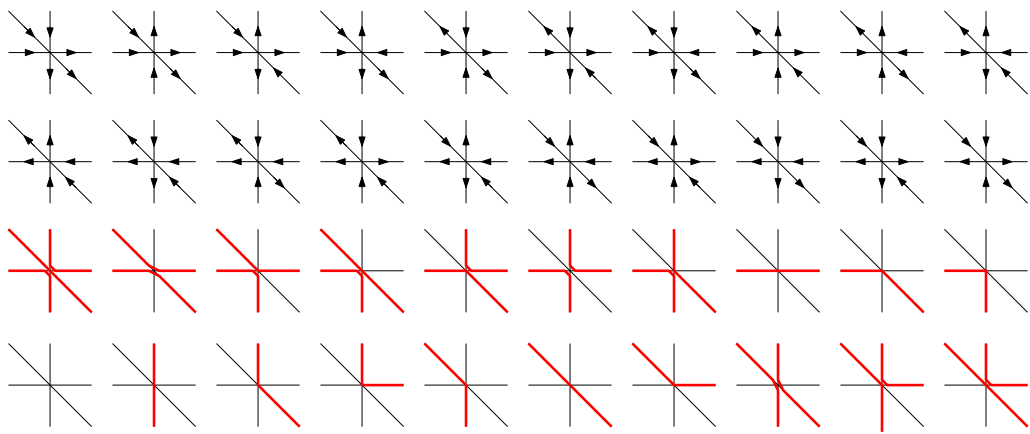}
\end{center}
\caption{\small The local vertex configurations of the 20V model (top) and their reformulation in terms of osculating Schr\"oder paths (bottom).}
\label{fig:twenty}
\end{figure}

Two-dimensional integrable lattice models of statistical physics have exhibited an extremely rich mathematical content, ranging from representation theory to probability and combinatorics. The present paper concentrates on the two-dimensional triangular lattice\footnote{Throughout the paper we view the two-dimensional triangular lattice in a sheared fashion: the vertices form a square lattice, and edges are those of the square lattice plus the second diagonal within each square. This allows for easier drawings, and is more adapted to our choices of domains.} 
version of the ice-type models, in the form of the so-called Twenty Vertex (20V) model \cite{Kel,Baxter}. The model is defined on a finite domain of the triangular lattice, by considering all possible orientations of edges of the elementary triangles, obeying the {\it ice rule}: ``there are equal numbers of 
incoming and outgoing edges at each vertex of the domain". This gives rise to the ${6\choose 3}=20$ local vertex configurations displayed 
in the top two rows of Fig.~\ref{fig:twenty}. The 20V model is really the triangular lattice version of the celebrated Six Vertex (6V) model \cite{Lieb} describing ice on a two-dimensional square lattice. 
The latter was at the core of the saga of Alternating Sign Matrices (ASM) \cite{Bressoud}, a remarkable sequence of connections between purely combinatorial objects such as Descending Plane Partitions (DPP) \cite{AndrewsDPP} in bijection with cyclically symmetric rhombus tilings of an hexagon with a central triangular hole \cite{KrattDPP}, Totally Symmetric Self-Complementary Plane Partitions (TSSCPP) \cite{tsscpp,zeilberger1} and two-dimensional models of statistical physics such as the 6V model \cite{IKdet,kuperberg1996another,tsu,kuperberg2002symmetry}, the O(1) loop model \cite{MNdG} and finally the Fully-Packed Loop model (FPL), at the center of the Razumov-Stroganov conjecture \cite{razustro}, later proved by Cantini and Sportiello \cite{CanSpo}. A striking feature of all the correspondences is the apparent absence of natural bijections between different classes of objects with the same cardinality $A_n=\prod_{j=0}^{n-1}\frac{(3j+1)!}{(n+j)!}=1,2,7,42,429,...$ The cases studied in this paper can be seen as other remarkable examples of coincidences of cardinalities of sets, between which no natural bijection is known.

\begin{figure}
\begin{center}
\includegraphics[width=12cm]{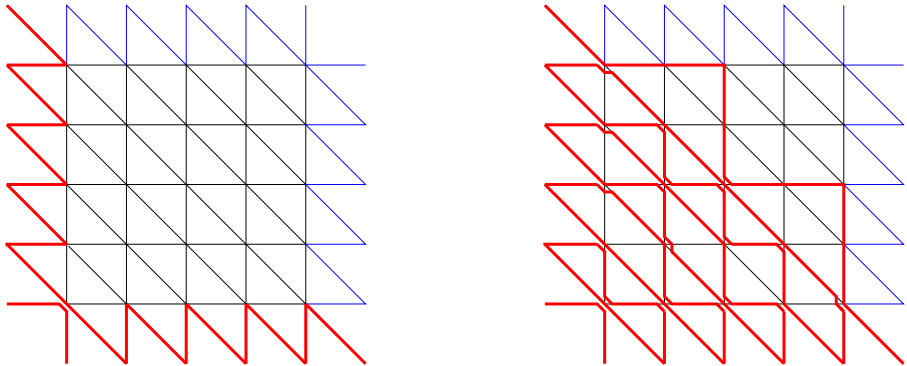}
\end{center}
\caption{\small DWBC1 boundary conditions for the 20V model and a sample configuration in the osculating Schr\"oder path formulation.}
\label{fig:dwbc1}
\end{figure}

In its physics literature debut, the 20V model was shown to be integrable provided some suitable choice of Boltzmann weights is made \cite{Kel,Baxter}, governed by 3 independent parameters, for which the thermodynamic free energy was computed. This used an explicit connection between the integrable 20V model on the triangular lattice and the integrable 6V model on the square lattice. 
Recently the model was revisited from a combinatorial perspective \cite{DG19}, by considering domains with so-called ``domain-wall boundary conditions" (DWBC), i.e. by imposing orientations
on the boundaries of the domain that force the existence of non-local fault lines within the domain, across which orientations are reversed (domain walls). In \cite{DG19}, several possible DWBC were considered (labeled 1,2,3,4) on a square $n\times n$ grid. Each type of DWBC gives rise to generalizations of ASM, coined Alternating Phase Matrices with entries among $0$ and the sixth roots of unity and conservation conditions along rows, columns and diagonal lines. 

The 20V model can be reformulated in terms of lattice paths, by picking a preferred orientation of all edges of the lattice 
(say right, down and diagonal down), and deciding that in any given orientation configuration, only edges with the preferred orientation are path steps. 
There are 3 possibilities for these steps: $(1,0),(0,-1),(1,-1)$ namely right, down and diagonal down. They form non-intersecting
but possibly ``kissing" paths (called osculating Schr\"oder paths) that are allowed to share vertices, where up to three paths bounce against each-other.
The 20 local configurations of osculating Schr\"oder paths are displayed in the two bottom rows of Fig.~\ref{fig:twenty}.

In the case of DWBC1,2 illustrated in Fig.~\ref{fig:dwbc1}\footnote{Ref. \cite{DG19} considers two choices labelled 1,2 for this first DWBC, which turn out to be equivalent modulo a $180^\circ$ rotation, and we show here a DWBC1 example.}, it was shown in \cite{DG19} that the partition function of the corresponding 20V model with uniform weights counts also the number of $2\times 1$ cyclically symmetric Domino Tilings of a Holey Aztec Square (HASDT), i.e. tilings of an Aztec-like square shaped domain of the square lattice with a central cross-shaped hole and with quarter-turn symmetry: this correspondence can be viewed as a generalization of the ASM-DPP correspondence \cite{BDZJ12}. The proof relies on the integrability of the 20V model allowing to transform the structure of the underlying lattice by moving lines around, and eventually relating the 20V partition function to that of the 6V model on a square grid with DWBC. As opposed to the case of ASM, which correspond to the 6V model with quantum parameter $q=e^{2i\frac{\pi}{3}}$, the 20V DWBC1,2 APM are enumerated by the partition function of a 6V DWBC model with $q=e^{i\frac{\pi}{8}}$. The connection between 20V and 6V being still valid in the presence of some non-trivial spectral parameters,
the result was generalized to include refinements as well.

\begin{figure}
\begin{center}
\includegraphics[width=12cm]{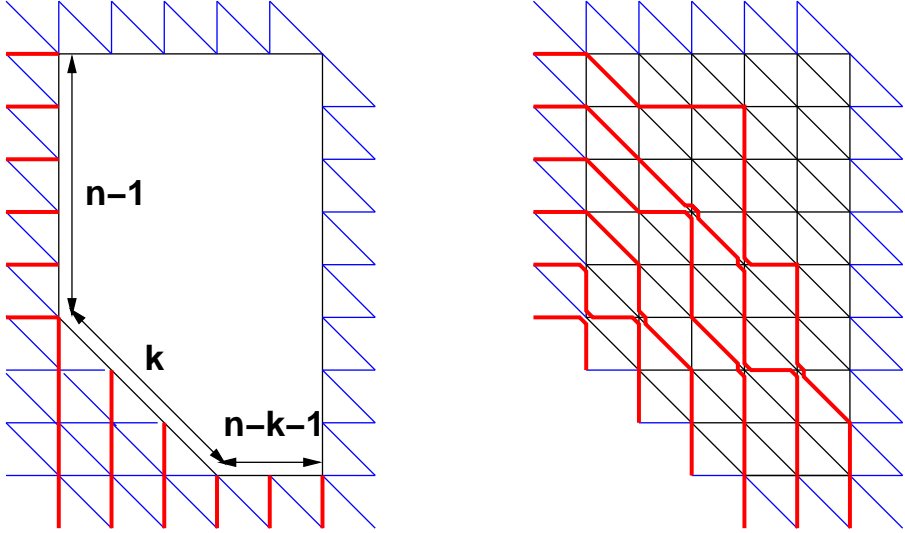}
\end{center}
\caption{\small Extended DWBC3 boundary conditions on the pentagon ${\mathcal P}_{n,k}$ for the 20V model and a sample configuration in the osculating Schr\"oder path formulation, with trivial forced vertices removed.}
\label{fig:extdwbc3}
\end{figure}

\begin{figure}
\begin{center}
\includegraphics[width=14cm]{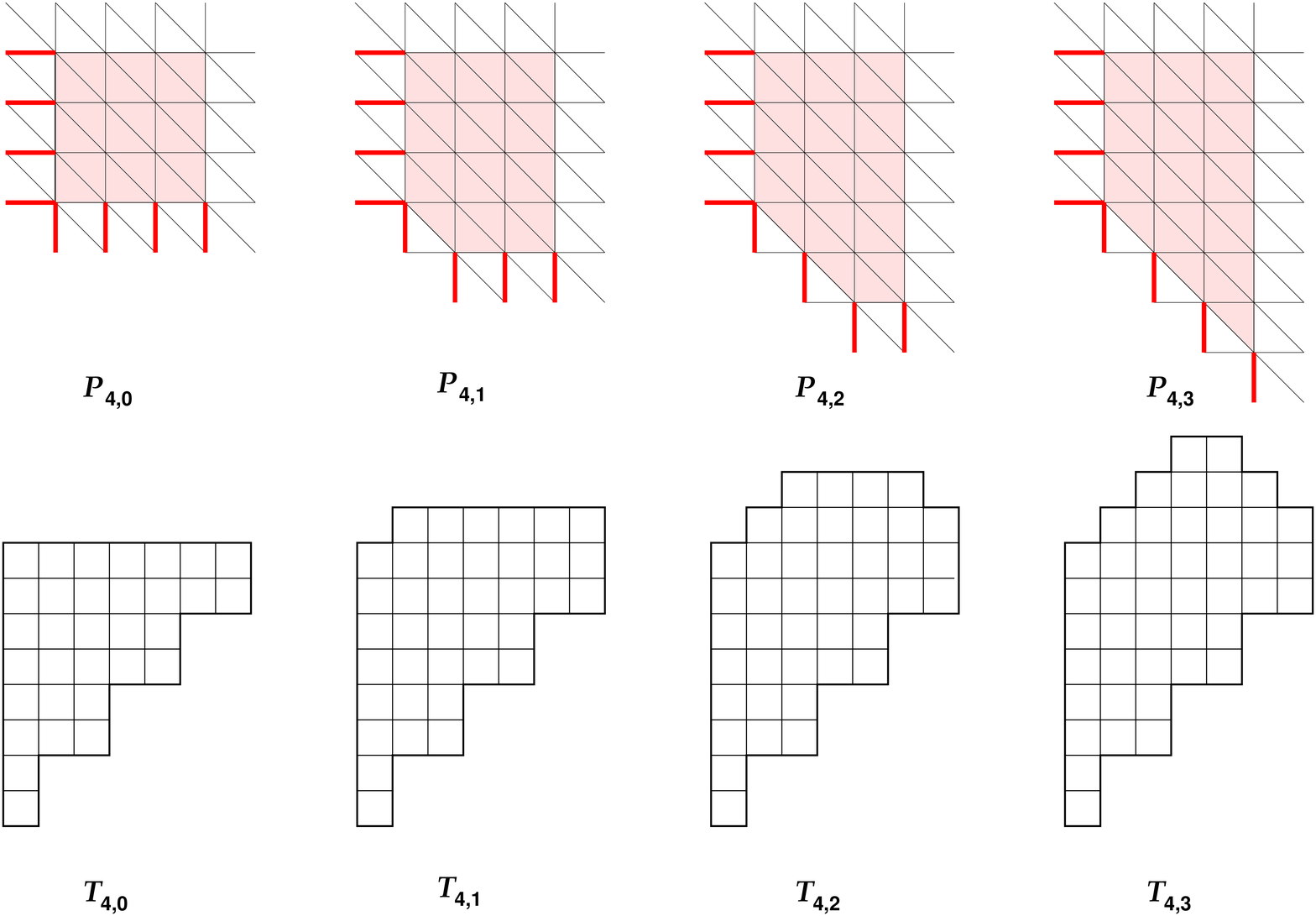}
\end{center}
\caption{\small Extended Pachter's triangular domains ${\mathcal T}_{n,k}$ for $n=4$ (bottom row), and the
corresponding $20V$ pentagonal domains ${\mathcal P}_{n,k}$ (top row).
The domain ${\mathcal T}_{n,k}$ must be
tiled by means of $2\times 1$ dominos, whereas ${\mathcal P}_{n,k}$ must carry configurations of the $20V$ 
model with the specified boundaries.}
\label{fig:pachtri}
\end{figure}

In the same paper, other combinatorial conjectures were made regarding the 20V DWBC3 model on a more general pentagonal grid
${\mathcal P}_{n,k}$, $0\leq k \leq n-1$, illustrated in Fig.~\ref{fig:extdwbc3}. It was conjectured in \cite{DG19} that the 20V DWBC3 configurations on ${\mathcal P}_{n,k}$ are equinumerous to the $2\times 1$ domino tiling configurations of extensions $\cT_{n,k}$ of the triangular domain $\cT_{n,0}$ of the square lattice introduced by Pachter \cite{Pachter} (see Fig.~\ref{fig:pachtri} for an illustration). In particular, for $k=0$ and $n=1,2,...$ the number of configurations of the 20V DWBC3 model on the square 
grid ${\mathcal P}_{n,0}$
forms the sequence $1,3,29,901,89893...$ which matches conjecturally that of domino tilings of the Pachter triangle $\cT_{n,0}$.

The aim of the present paper is to shed some light on these latter conjectures and to prove some of them. More precisely, we will concentrate on the 20V DWBC3 model on the maximal domains $\cQ_n:=\cP_{n,n-1}$ (as in Fig.~\ref{fig:pachtri} top right), whose configuration numbers give rise to the sequence
$1,4,60,3328,678912...$ and prove that they are counted by the numbers of $2\times 1$ domino tilings of the maximally extended Aztec triangle $\cT_n:=\cT_{n,n-1}$ (as in  Fig.~\ref{fig:pachtri} bottom right) . 

\subsection{Outline of the paper and main results}

The paper is organized as follows.

The partition function of the 20V model is computed in Section \ref{pfsec}. Using integrability, we transform the 20V model into a 6V model on a rectangular grid with different boundary conditions, the so-called U-turn boundary conditions\footnote{Remarkably, the same model but with the quantum parameter $q=e^{2i\frac{\pi}{3}}$ was shown to enumerate Vertically Symmetric ASM (VSASM) \cite{kuperberg2002symmetry}.} \cite{kuperberg2002symmetry}, and with quantum parameter $q=e^{i\frac{\pi}{8}}$. The fully inhomogeneous version of the 20V model involves a large number of arbitrary ``spectral" parameters, which can be specialized to the {\it combinatorial point} where all the Boltzmann weights are $1$. On the other hand, the 6V U-turn boundary condition partition function is expressed as an explicit determinant involving the spectral parameters of the model (Section \ref{kupsec}).
The main difficulty in computing the 20V model partition function is to derive the homogeneous limit of this determinant, when all spectral parameters tend to their combinatorial point values. The result takes the form of a compact determinant expression for the total number $Z_n^{20V}$ of configurations of the 20V model on the domain $\cQ_n$ (Theorem \ref{thm20V}):

\begin{thm} The total number $Z_n^{20V}$ of configurations of the 20V DWBC3 model on the domain $\cQ_n$ reads:
$$
Z_n^{20V}=\det_{0\leq i,j\leq n-1}\left(\left.
\frac{(1+u^2)(1+2u-u^2)}{(1-u^2v)\Big((1-u)^2-v(1+u)^2\Big)}\right\vert_{u^iv^j} \right)
$$
\end{thm}

The enumeration of domino tilings $Z_n^{DT}$ of the Aztec triangle $\cT_n$ is performed in Section \ref{dtsec} by a bijective formulation in terms of non-intersecting Schr\"oder lattice paths, leading to a compact determinant formula (Theorem \ref{thmAT}): 

\begin{thm} The total number $Z^{DT}_n$ of domino tilings of the Aztec triangle $\cT_n$ reads:
$$
Z^{DT}_n=\det_{0\leq i,j \leq n-1}\left( \left.
\frac{1+u}{1-v -4 u v - u^2 v +u^2 v^2}\right\vert_{u^iv^j} \right)$$
\end{thm}

The identity between 20V and domino tiling enumerations is proved in Section \ref{correspsec} by noting that both enumeration results of Theorems \ref{thm20V} and \ref{thmAT} express the numbers as the determinant of a finite truncation of infinite matrices $P$ and $M$ respectively, and by exhibiting an explicit infinite lower triangular matrix $L$ with entries 1 on the diagonal such that $L\,M=P$ (Theorem \ref{equivalencethm}):

\begin{thm} The partition function of the 20V DWBC3 model on the quadrangle ${\mathcal Q}_n$ and that of the domino tilings of the Aztec triangle 
$\cT_n$ coincide for all $n\geq 1$.
\end{thm}

The second part of the paper is devoted to a refinement of the enumeration result, in which we keep track of one particular statistic on the set of configurations. For 20V configurations, we record the position of the first visit to the rightmost vertical column. For domino tilings, we record the first visit to the line of slope $-1$ through the topmost endpoint of the paths. The relation between the two refinements is obtained by repeating the analysis of the previous sections. 

In Section \ref{ref20Vsec}, we follow the transformation of the 20V model into the U-turn 6V model while keeping an arbitrary spectral parameter $w$ in the last column, all other parameters being sent to their combinatorial point values. This gives rise to a compact determinant formula for the generating polynomial of refined 20V partition functions in the variable $\tau=\frac{q^{-2}-q^2 w}{q^{2}-q^{-2} w}$ (Theorems \ref{oneptthm} and \ref{refcor}).  

In Section \ref{refdtsec}, we use again the non-intersecting Schr\"oder lattice path formulation of the domino tiling problem to derive a compact determinant formula for a generating polynomial of the refined partition functions (Theorem \ref{refdominothm}). 

In Section \ref{seclutfin}, the results of Theorems \ref{oneptthm} and \ref{refdominothm} are shown to be equivalent (Theorem \ref{fullequivthm}). The proof, inspired by techniques employed in \cite{BDZJ12,DFLAP} relies on the observation that both results are {\it perturbations} of the uniform case, in which both infinite matrices $M,P$ are modified into infinite matrices ${\bar M}^{(n)}(\tau),P^{(n)}(\tau)$ by an additional term affecting only their $n$-th and higher columns, and for which the relation $L\, {\bar M}^{(n)}(\tau)=P^{(n)}(\tau)$ still holds.

Section \ref{conjsec} is devoted to an investigation of the other conjectures regarding the 20V model on $\cP_{n,n-k}$ and the domino tilings of $\cT_{n,n-k}$ for $k>1$. We present direct proofs of the conjectures for $k=2$ and $k=3$ based on recursion relations involving the $k=1$ result and its refinements.

A few concluding remarks are gathered in Section \ref{concsec}. In Section \ref{secconj}, we present the following:
\begin{conj}\label{numconj0}
The total numbers $Z^{20V}_n=Z_n^{DT}$ of configurations of the 20V DWBC3 model on the quadrangle $\cQ_n$ and of domino tilings of the Aztec triangle $\cT_n$ read: 
\begin{equation}\label{conjnum0}
Z^{20V}_n=Z_n^{DT}= 2^{n(n-1)/2}\prod_{i=0}^{n-1}  \frac{(4i+2)!}{(n+2i+1)!}
\end{equation}
\end{conj}
\noindent for which we have no proof at this time. In Section \ref{altsec}, we present alternative formulas for the numbers $Z^{20V}_n=Z_n^{DT}$, including a determinant involving binomials (Theorem \ref{binothm}) and a constant term formula (Theorem \ref{ctthm}), in the same spirit as the TSSCP enumeration of Zeilberger \cite{zeilberger1}. In Section \ref{stepweightsec}, we show how a simple modification of the results of this paper allows to include an extra weight $\gamma$ per horizontal step in the non-intersecting Schr\"oder lattice path formulation of the domino tiling problem. Section \ref{seconc} is the conclusion {\it per se}.

\bigskip

\noindent{\bf Acknowledgments.} We acknowledge partial support from the Morris and Gertrude Fine endowment and the NSF grant DMS18-02044.



\section{Partition function of the 20V model}
\label{pfsec}

\subsection{The 20V model: integrable weights}

Given a domain $D$ of the two-dimensional triangular lattice, we consider the set of orientations of all lattice edges within $D$ that satisfy the ice rule, 
namely that {\it each vertex is incident to exactly 3 incoming and 3 outgoing edges}. Boundary edge orientations are fixed, as part of the definition of the model. 
As mentioned before, the 20V model configurations are alternatively described by osculating Schr\"oder paths on edges of the domain $D$ that travel right and down. The boundary edges of $D$ are fixed to be occupied/empty as part of the definition of the model. The occupied boundary edges split into equal numbers of starting and endpoints of the paths.

\begin{figure}
\begin{center}
\includegraphics[width=12cm]{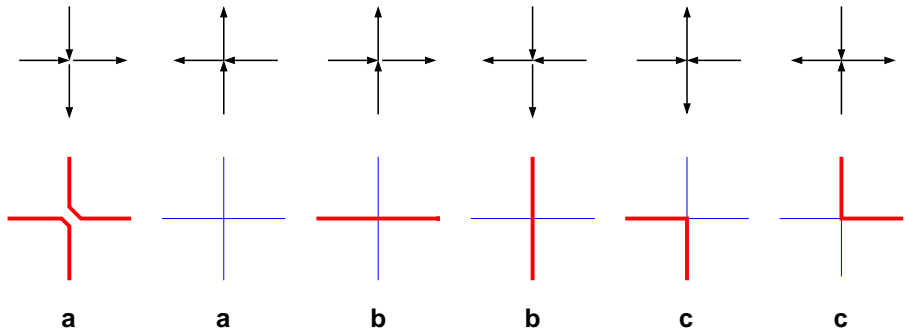}
\end{center}
\caption{\small The local vertex configurations of the 6V model (top row) and the corresponding osculating path configurations (bottom row). We have indicated the weights $a,b,c$ of the respective configurations.}
\label{fig:six}
\end{figure}

Let us now describe the weighting of the configurations. Each vertex receives a local weight according to its configuration, and each configuration on $D$ receives a weight equal to the product of the local weights of its inner vertices. The weights are chosen in a special manner, which we describe now.
We attach complex ``spectral parameters" to the lines of the lattice that intersect $D$. If a vertex is at the intersection of three such lines, say with parameters $z,t,w$, we attach a weight $W_i(z,t,w)$, $i=1,2...,20$ according to the list in Fig.~\ref{fig:twenty}.
In \cite{Kel,Baxter,DG19}, the weights are built out of elementary pieces as follows. 

Consider a triple intersection between a horizontal, vertical and diagonal line, and translate slightly the diagonal, say to the right: $\raisebox{-.5cm}{\hbox{\epsfxsize= 4.cm \epsfbox{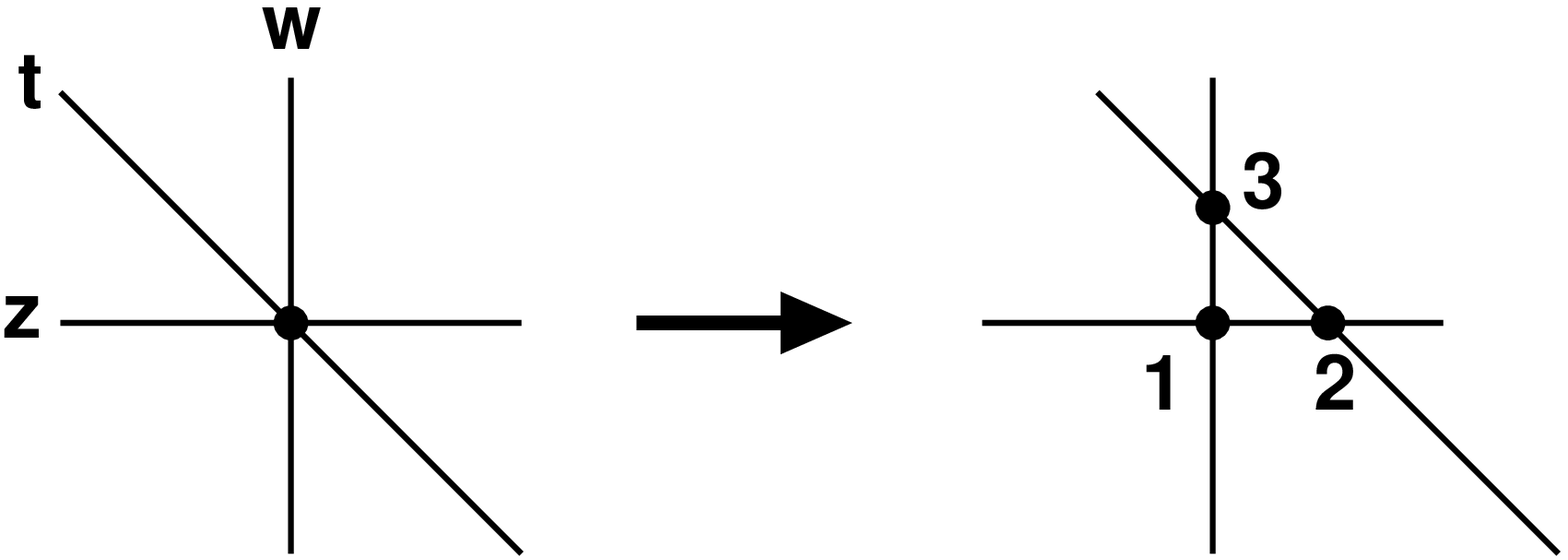}}}$.
This amounts to resolving the triple intersection into three simple intersections, labeled $1,2,3$. 
\begin{figure}
\begin{center}
\includegraphics[width=8cm]{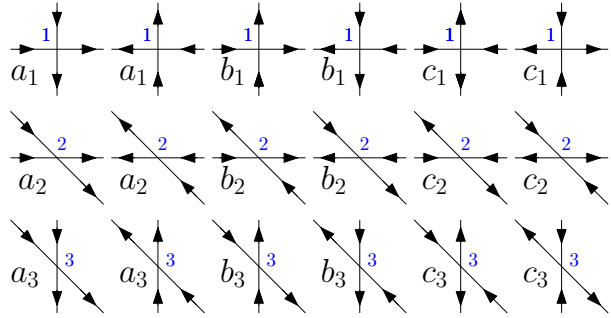}
\end{center}
\caption{\small The local vertex configurations of the three 6V models after resolving all triple intersections in the 20V model.}
\label{fig:kago}
\end{figure}
Assume the edges were originally oriented and obeying the ice rule. 
The resolution above keeps the six outer edge orientations and creates 3 new inner edges,
which may be oriented in such a way that the ice rule is now satisfied at each simple intersection, which now has 4 adjacent edges. At each simple intersection, this gives rise to ${4\choose 2}=6$ possible local configurations, those of the so-called six vertex (6V) model (see Fig.~\ref{fig:six} top for a list of the vertex environments and their 3 customary weights $a,b,c$ that are invariant under reversal of all orientations). Here we are really dealing with three different 6V models, attached to the three simple intersections. The 6V weights are parameterized by the two spectral parameters of their two lines, and an additional index 1,2,3 corresponding to respectively a horizontal-vertical, horizontal-diagonal, and diagonal-vertical intersection. 
The three sets of 6V weights are denoted by $a_i,b_i,c_i$ and are functions of their two spectral parameters 
(see Fig.~\ref{fig:kago} for the list of all 18 possible configurations and their weights).
The 20V weights can now be defined in terms of the 6V weights of the resolution as follows:
$$W^{20V}(z,t,w):= \sum_{{\rm inner}\  {\rm edge}\ {\rm configs}.} w^{6V}_1(z,w) \, w^{6V}_2(z,t)\, w^{6V}_3(t,w) $$
where the sum extends over all allowed inner edge orientations, and $w^{6V}_i\in \{a_i,b_i,c_i\}$ stands for the 6V weight of the corresponding local  configuration. For instance, the resolution:
\begin{equation}\label{casex} \raisebox{-.5cm}{\hbox{\epsfxsize= 9.cm \epsfbox{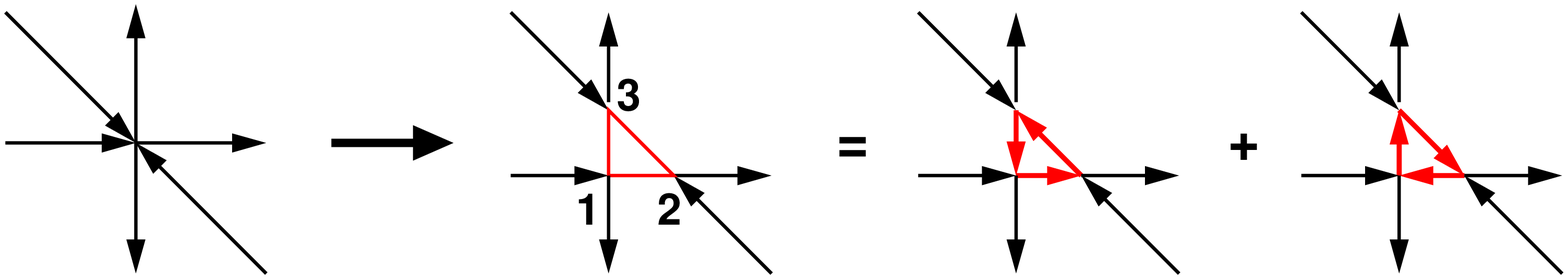}}} \end{equation}
corresponds a 20V weight $\omega=a_1b_2c_3+c_1c_2b_3$.

The weights are further constrained by the integrability condition which amounts to requiring that these values be independent of the resolution. For instance, we could have moved the diagonal slightly to the left instead of the right, giving rise to an {\it a priori} different definition of the $W^{20V}$s. 
The requirement that the two definitions coincide results in a system of cubic relations for the 6V weights (the celebrated Yang-Baxter equations). In the case \eqref{casex}, the Yang-Baxter equation reads
$\omega=a_1b_2c_3+c_1c_2b_3=b_1a_2c_3$. We use the solution of these equations from \cite{DG19}, equivalent to that of \cite{Kel,Baxter}.
Define the following three functions of the two spectral parameters $z,w$:
\begin{equation} \label{w6v}
A(z,w)=z-w,\quad B(z,w)=q^{-2}z -q^{2} w,\quad C(z,w)=(q^2-q^{-2}) \sqrt{z w} ,
\end{equation}
and pick the following 6V weights:
\begin{equation} \label{3w6v}
\small{
\begin{matrix}
a_1=A(z,w)=z-w, \hfill &b_1=B(z,w)=q^{-2}z -q^2 w,\hfill &c_1=C(z,w)=(q^2-q^{-2}) \sqrt{z w} \hfill\\
a_2=A(qz,q^{-1}t)=q z-q^{-1}t, \hfill &b_2=B(qz,q^{-1}t)=q^{-1}z -q t,\hfill &c_2=C(qz,q^{-1}t)=(q^2-q^{-2}) \sqrt{z t} \hfill\\
a_3=A(qt,q^{-1}w)=q t-q^{-1}w, \hfill &b_3=B(qt,q^{-1}w)=q^{-1}t -q w,\hfill &c_3=C(qt,q^{-1}w)=(q^2-q^{-2}) \sqrt{t w}\hfill
\end{matrix}}
\end{equation}

\begin{figure}
\begin{center}
\includegraphics[width=15cm]{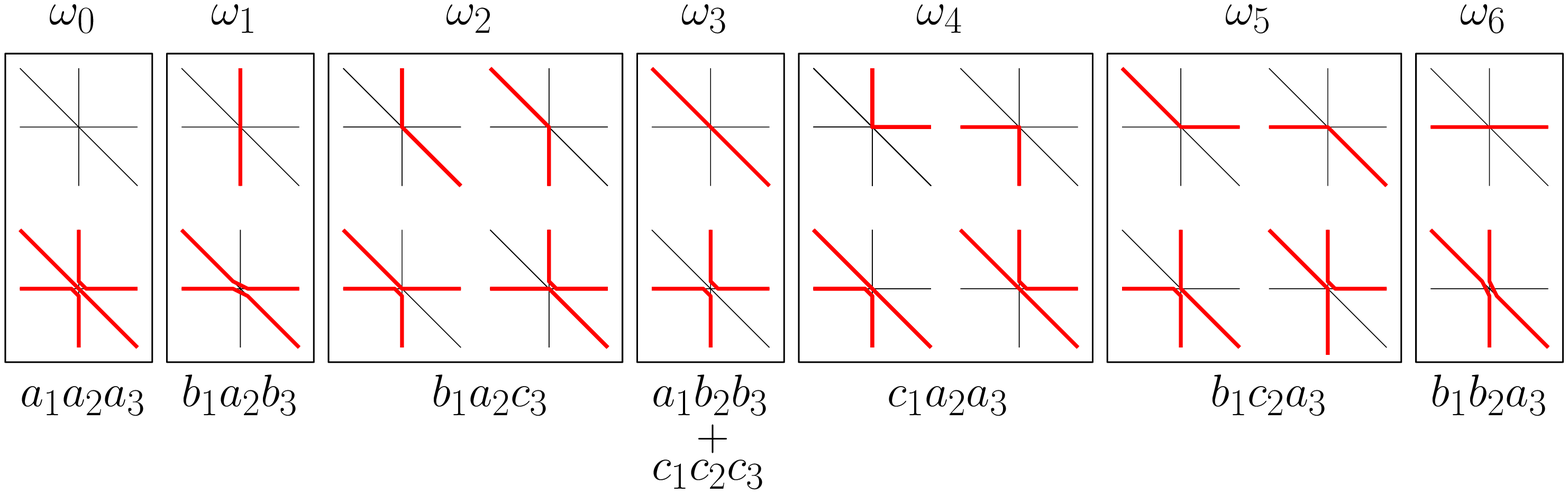}
\end{center}
\caption{\small Integrable weights of the 20V model and their expression in terms of weights of the three sublattice 6V models.}
\label{fig:20vweights}
\end{figure}

It is easy to show that the Yang-Baxter equations are satisfied by these, and lead to the following expressions for the 20V integrable weights:
\begin{eqnarray}\label{weights20V}
\omega_0&=& (z-w)(qz-q^{-1}t)(qt-q^{-1}w)\nonumber \\
\omega_1&=&(q^{-2}z-q^2w)(qz-q^{-1} t)(q^{-1}t-qw) \nonumber \\
\omega_2&=&(q^{-2}z-q^2w)(qz-q^{-1}t)(q^2-q^{-2})\sqrt{tw}\nonumber \\
\omega_3&=& z t w (q^2-q^{-2})^3+(z-w)(q^{-1}z -q t)(q^{-1}t-q w) \nonumber \\
\omega_4&=&(q^2-q^{-2})\sqrt{zw}(qz-q^{-1}t)(qt-q^{-1}w)\nonumber \\
\omega_5&=&(q^{-2}z-q^2w)(q^2-q^{-2})\sqrt{zt}(qt-q^{-1}w)\nonumber \\
\omega_6&=&(q^{-2}z-q^2w)(q^{-1}z-q t)(qt-q^{-1}w)
\end{eqnarray}
where the index refers to the vertex configurations of Fig.~\ref{fig:20vweights} (our running example \eqref{casex} corresponds to $\omega=\omega_2$).

We will be first interested in the pure enumeration of the configurations of the 20V model on specific domains with prescribed boundary conditions, which requires uniform weights for all vertices. As it turns out, there exists a choice of the spectral parameters and of the parameter $q$ which ensures that all the weights \eqref{weights20V} are equal to $1$.

\begin{defn}
The ``combinatorial point" where all $\omega_i=1$ corresponds to the values 
\begin{equation}\label{combipoint}
q=e^{i\pi/8} \quad {\rm and}\quad (z,t,w)=\alpha (q^6,1,q^{-6}), \quad \alpha =2^{-5/6} q^{-4} 
\end{equation}
\end{defn}

The integrability property of the weights \eqref{weights20V} will be crucially used in this paper to transform partition functions. The property can be used to deform/move any line of the lattice without changing the value of quantities such as the partition function, i.e. the weighted sum over configurations. For that, the presence of arbitrary spectral parameters is essential. Combinatorial (enumeration) results will be obtained by eventually sending all parameters to their combinatorial point values \eqref{combipoint}.

\subsection{Boundary conditions}

\begin{figure}
\begin{center}
\includegraphics[width=6cm]{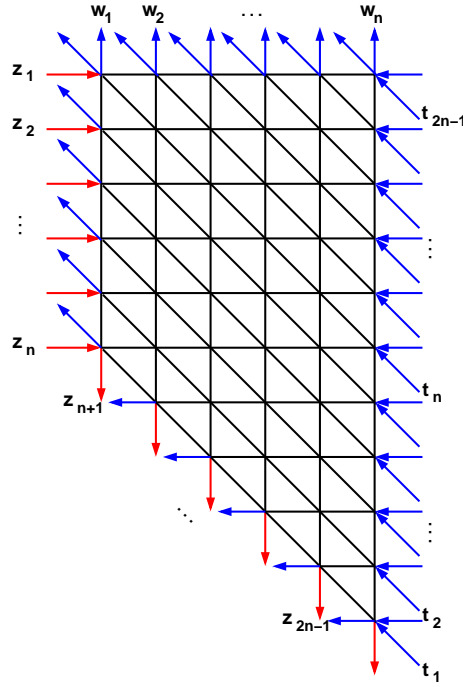}
\end{center}
\caption{\small Quadrangular domain ${\mathcal Q}_n$ for the partition function of the 20V model. We have indicated the horizontal spectral parameters $z_1,z_2,...,z_{2n-1}$, the vertical ones $w_1,w_2,...,w_n$ and the diagonal ones $t_1,t_2,...,t_{2n-1}$. We impose the following boundary conditions: arrows pointing towards the domain on the E boundary, out of the domain on the N and SW boundaries, and horizontal arrows point inward and diagonal arrows point outward along the W boundary.}
\label{fig:domain}
\end{figure}

We now consider the configurations of the 20V model on the quadrangular domain ${\mathcal Q}_n={\mathcal P}_{n,n-1}$ depicted in Fig.~\ref{fig:domain}, with a horizontal north boundary (N), vertical west and east (W,E) boundaries and diagonal south-west (SW) boundaries.
We assign spectral parameters $\bz=(z_1,z_2,...z_{2n-1})$ to horizontal lines (top to bottom), $\bw=(w_1,w_2,...,w_n)$
to vertical lines (left to right) and $\bt=(t_1,t_2,...,t_{2n-1})$ to diagonal lines (bottom to top). 
The corresponding fully inhomogeneous partition function is denoted by $Z_n^{20V}[\bz;\bt;\bw]$: it is the sum over all configurations of the 20V model compatible with the imposed boundary conditions, of the product of corresponding local weights 
$\omega_m(z_i,t_j,w_k)$ at the intersection of the three horizontal, diagonal and vertical lines carrying spectral 
parameters $z_i,t_j,w_k$ respectively.

We also consider a semi-homogeneous version of this partition function
$Z_n^{20V}[z,t;\bw]$ obtained by setting $z_i=z$ and $t_i=t$ for all $i$, while keeping the $w_i$ arbitrary.

The total number $Z_n^{20V}$ of configurations is obtained by considering the abovementioned combinatorial point, 
namely further restricting $w_i=w$ for all $i$, and using the values \eqref{combipoint}.

\subsection{Transformation to a mixed 6V model}
\label{transfosec}

\begin{figure}
\begin{center}
\includegraphics[width=15cm]{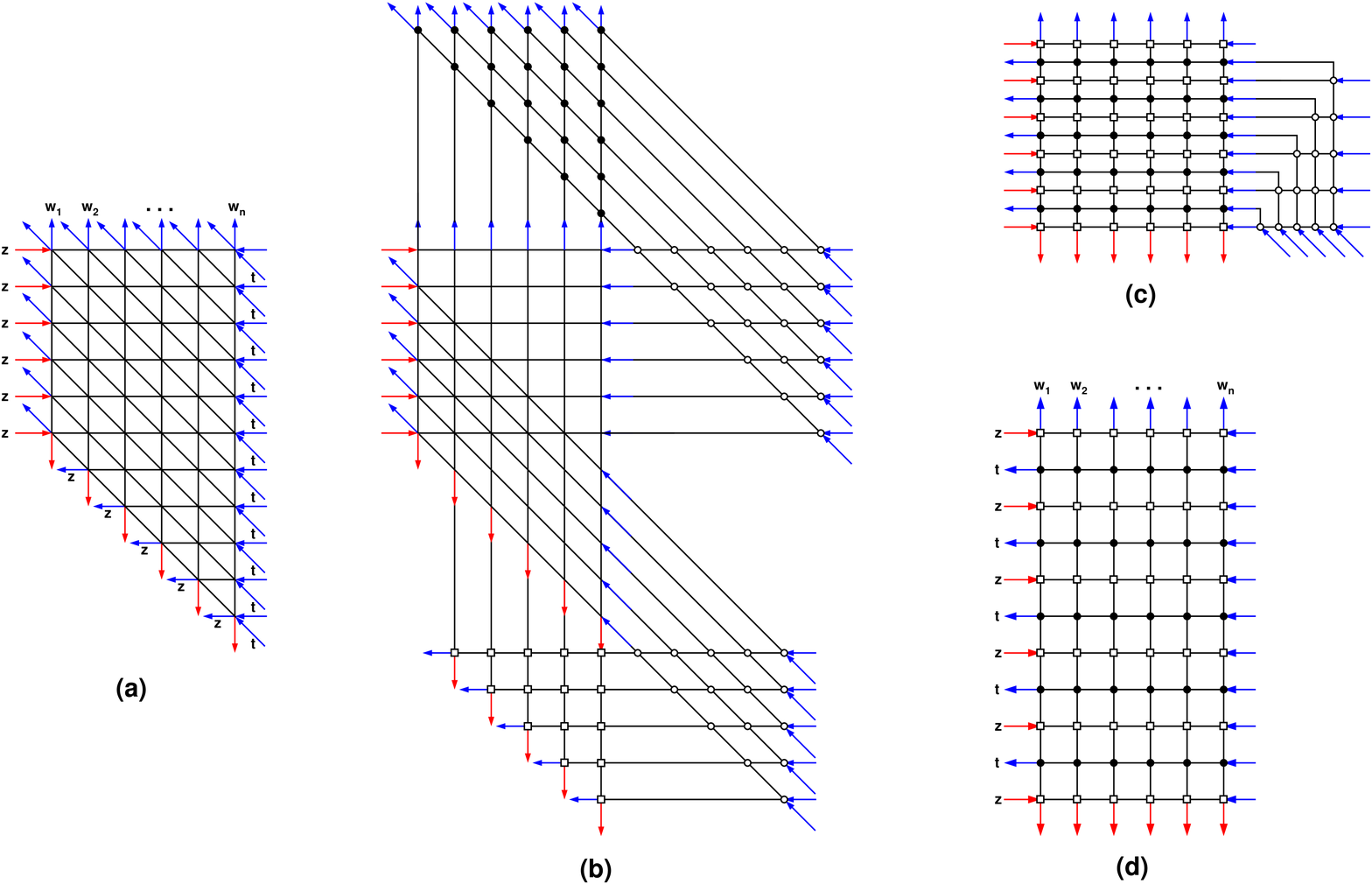}
\end{center}
\caption{\small Transformation of the original 20V model partition function $Z_n^{20V}[z,t;\bw]$ (a) into 
the alternating boundary mixed 6V model partition function $Z_{n}^{6V}[z,t;\bw]$ (d) (represented here for $n=6$), 
by moving the diagonal lines, using the Yang-Baxter equation. The marking of vertices 
corresponds to the weights of the three 6V models on the three sublattices 1 (empty square), 2 (empty circle) and 3 (filled circle).
In (a), spectral parameters have uniform values $z$ (horizontal), $t$ (diagonal), and arbitrary values $w_i$ (vertical). As a consequence, in (d)  the horizontal spectral parameters alternate between $z,t$ and the vertical ones are $w_1,w_2,...,w_n$.}
\label{fig:transfo}

\end{figure}
The semi-homogeneous partition function $Z_n^{20V}[z,t;\bw]$
can be identified with that of a special mixed 6V model on a rectangular grid,
by use of the Yang-Baxter equation, which allows to displace spectral lines while keeping track of their intersection 
pattern (see Fig.~\ref{fig:transfo} (a-d)). 

Starting from (a), we first move up the $n$ uppermost diagonal lines, and down the $n-1$
lowermost horizontal lines. Keeping track of their intersections with the rest of the grid lines, these give rise to four triangular 
domains with {\it fixed} arrow configurations (all with identical arrow orientations along lines), due to the propagation of the boundary conditions via the ice rule. The two top domains have $n(n+1)/2$ vertices each, while the two bottom ones have $n(n-1)/2$.

Erasing them all, we get a new partition function. The original one is the product of the trivial weights of all the erased vertices,
namely: $ a_2(z,t)^{n^2} \prod_{i=1}^n b_1(z,w_i)^{i-1} a_3(t,w_i)^i$ times the new partition function with all the trivial vertices erased.
(Here $a_i,b_i,c_i$ denote the uniform weights of the three sublattices $i=1,2,3$ 6V models \eqref{3w6v}).
The latter partition function is further transformed into (c) by straightening the remaining diagonal lines, and keeping track of their intersections.
This gives rise to a last triangular domain (on the right), again with fixed trivial arrow configurations. Erasing the corresponding
$n(n-1)/2$ trivial vertices gives rise to the partition function of (d), which we denote by $Z_{n}^{m6V}[z,t;\bw]$, 
up to an extra factor
of $a_2(z,t)^{n(n-1)/2}$. The acronym $m6V$ stands for mixed 6V model partition function: it is defined on a rectangular grid of size $(2n-1)\times n$ with 6V weights of sublattices 1 and 3 alternating on horizontal lines,
respectively with uniform horizontal spectral parameters $z$ and $t$, and vertical spectral parameters $w_1,w_2,...,w_n$, and with arrows pointing alternatively in and out along the W border.

Recording all the weights of the erased vertices, we get the relation:

\begin{equation}\label{relamod}
Z_n^{20V}[z,t;\bw]= \left(a_2(z,t)^{n(3n-1)/2} \, \prod_{i=1}^n b_1(z,w_i)^{i-1} \, a_3(t,w_i)^{i}\right)\, Z_{n}^{m6V}[z,t;\bw]
\end{equation}


\subsection{Transformation to a U-turn boundary 6V model}
\label{usec}

\begin{figure}
\begin{center}
\includegraphics[width=13cm]{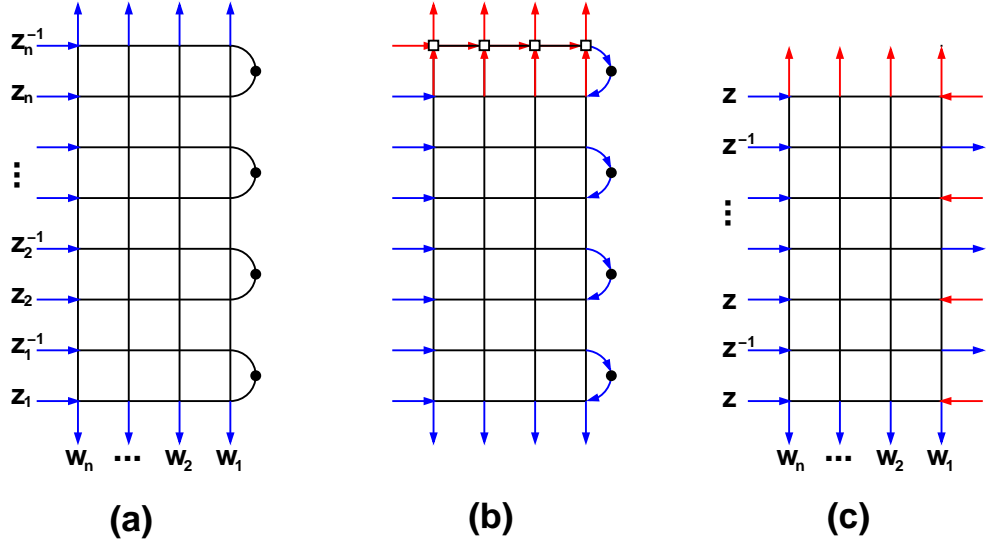}
\end{center}
\caption{\small The U-turn boundary 6V model inhomogeneous partition function (a) reduces when all  $z_i=z=p^{-2}q^{-2}$  to the configurations
of (b) where all the arrow orientations at the U-turns point down, and where the top vertices are all frozen in a b-type configuration. Erasing the latter and cutting the U-turns leads to the partition function (c), in which horizontal spectral parameters alternate between $z$ and $z^{-1}$ from bottom to top, while the vertical ones are $w_n,...,w_2,w_1$ from left to right.
}

\label{fig:Uboundary}

\end{figure}

The U-turn boundary 6V model was independently studied by Tsuchiya \cite{tsu} and Kuperberg \cite{kuperberg2002symmetry}. The latter uses the 6V weights
$(A,B,C)(qz,q^{-1} w)/\sqrt{z w}$ (and up to a renaming of variables $x=\sqrt{z/w}$, $a=q$, $b=p$). 
Here we drop the overall factor of $1/\sqrt{z w}$ in the definition of the weights, which amounts to using the 6V weights \eqref{3w6v}
of the sublattices 2 or 3.

Let 
$Z_{n}^{6V-U}[z_1,z_2,...,z_{n};w_1,w_2,...,w_n;p]\equiv Z_{n}^{6V-U}[\bz;\bw;p]$ 
denote the partition function of the inhomogeneous 6V model
on a rectangular grid of size $2n\times n$ with alternating horizontal spectral parameters $z_1,z_1^{-1},...,z_n,z_n^{-1}$ 
(from bottom to top), vertical spectral parameters $w_n,w_{n-1},...,w_1$ (from left to right), with entering arrows along the W boundary, exiting arrows along the N and S boundaries, and U-turn
boundaries on the E, as described in Fig.~\ref{fig:Uboundary} (a).
The marked dot on the U-turn indicates that the
spectral parameter changes from $z_i$ (bottom) to $z_i^{-1}$ (top) and that the configuration of the corresponding U-turn 
receives a weight $Y(z_i)$, where
\begin{equation}\label{Uweight}
Y(z)=\left\{ \begin{matrix} pq^{-1}-p^{-1}q z & \raisebox{-.3cm}{\hbox{\epsfxsize=.6cm \epsfbox{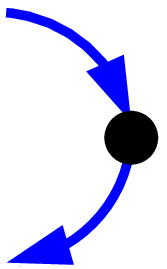}}}\\
\ & \  \\
pqz-p^{-1}q^{-1}& \raisebox{-.3cm}{\hbox{\epsfxsize=.6cm \epsfbox{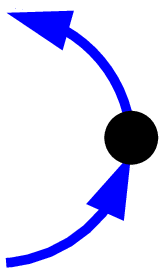}}}
\end{matrix}\right.
\end{equation}  
These weights obey Sklyanin's Reflection Yang-Baxter equation \cite{sklyanin}.
We now relate the partition function of this model to that of the mixed 6V model of previous section.

Setting all $z_i=z=p^{-2}q^{-2}$ forces each U-turn to have an arrow pointing down, as indicated in Fig.~\ref{fig:Uboundary} (b).
However this imposes no constraint on the vertical spectral parameters $w_i$, which we keep arbitrary.
In the resulting semi-homogeneous partition function, we
may then cut out all the marked dots on the U-turns, and moreover use the ice rule to dispose of the top-most 
horizontal line,
which gives rise to $n$ vertices with trivial b-type configurations. Erasing
these while keeping track of their weight leaves us with the semi-homogeneous partition function $Z_{n}^{6V}[z;\bw]$ corresponding to
Fig.~\ref{fig:Uboundary} (c), with alternating uniform horizontal spectral parameters $z,z^{-1}$, and arbitrary vertical spectral parameters $\bw=(w_1,w_2,...,w_n)$:
\begin{equation}
\label{homog}
Z_{n}^{6V}[z;\bw]:=  \frac{Z_{n}^{6V-U}[z;\bw;p]}{(p^2-p^{-2})^n\prod_{i=1}^n (p-p^{-1} w_i)}  \ , \quad p=z^{-1/2}q^{-1},
\end{equation}
where we have divided by the top row b-type weights $q^{-1}z^{-1}-q w_i=q p( p-p^{-1} w_i)$ and by the U-turn weights 
$pq^{-1}-p^{-1}q z=p^{-1}q^{-1}(p^2-p^{-2})$, using the relation $z=p^{-2}q^{-2}$.

The partition function $Z_{n}^{6V}[z;\bw]$ \eqref{homog} may now be identified with $Z_{n}^{m6V}[z,t;\bw]$ of the previous section as follows. 
Note first that a rotation of $180^\circ$
maps the configurations of Fig.~\ref{fig:Uboundary} (c) onto those of Fig.~\ref{fig:transfo} (d) for general size $n$, and then 
that this rotation flips all the local arrow orientations of the vertices, which therefore keep the same weights \eqref{w6v}.
However a discrepancy subsists between the weights in $Z_{n}^{m6V}[z,t;\bw]$ and those in $Z_{n}^{6V}[z;\bw]$.
Indeed, the 6V weights for $Z_{n}^{m6V}[z,t;\bw]$ alternate between those of the sublattice 1 on odd rows and 3 on even rows.
The weights on the even rows are those of the sublattice 3 in both models, and we must therefore identify their spectral parameters by setting $t=z^{-1}$.
On odd rows, using the relation $(a_1,b_1,c_1)(z,w_i)=(a_3,b_3,c_3)(q^{-1} z,q w_i)=q (a_3,b_3,c_3)(q^{-2} z,w_i)$, we may identify these weights with
$(a_3,b_3,c_3)(z,w_i)$ up to the overall factor $q$, provided we  redefine odd row parameters $z\to q^2 z$ (while keeping $t=z^{-1}$ on even rows). 
This results in the following relation between the semi-homogeneous partition functions:
$$ Z_{n}^{m6V}[q^2 z,z^{-1};\bw]= q^{n^2} Z_{n}^{6V}[z;\bw] $$
with $z=p^{-2}q^{-2}$.
Together with \eqref{relamod}, this leads to the relation:
\begin{equation}\label{twentytosix}
Z_n^{20V}[q^2z,z^{-1};\bw]= \frac{q^{n^2}a_2(q^2z,z^{-1})^{n(3n-1)/2} \, 
\prod_{i=1}^n b_1(q^2z,w_i)^{i-1} \, a_3(z^{-1},w_i)^{i}}{(p^2-p^{-2})^n\prod_{i=1}^n (p-p^{-1} w_i)}\, Z_{n}^{6V-U}[z;\bw;p]
\end{equation}
with $p=z^{-1/2}q^{-1}$.

To reach the combinatorial point \eqref{combipoint}, let us pick uniform 
20V spectral parameters so that $(q^2 z,z^{-1},w_i)=\delta (q^6,1,q^{-6})$, namely $z=q^2,\delta=q^{-2},p=q^{-2}$
and $w_i=q^{-8}$ for all $i$,
and finally set $q=e^{i\pi/8}$. Using the homogeneity of the weights, this allows to relate:
\begin{eqnarray}\label{20vone}
Z_n^{20V}&=&\left(\frac{\al}{\delta}\right)^{\frac{3n(3n-1)}{2}}\, Z_n^{20V}[q^2 z,z^{-1},w]\big\vert_{z=q^2;w=q^{-8}} \nonumber \\
&=& \frac{b_1(q^4,q^{-8})^{\frac{n(n-1)}{2}}a_2(q^4,q^{-2})^{\frac{n(3n-1)}{2}}  a_3(q^{-2},q^{-8})^{\frac{n(n+1)}{2}}q^{n^2}}{\left(2^{5/6}q^{2}\right)^{\frac{3n(3n-1)}{2}}(q^{-4}-q^4)^n (q^{-2}-q^{-6})^n} \,Z_{n}^{6V-U}[q^2,q^{-8};q^{-2}]\nonumber \\
&=&  2^{-\frac{n(5n+3)}{4}}\, q^{-5n} \,Z_{n}^{6V-U}[q^2,-1;q^{-2}]
\end{eqnarray}
with $\al$ as in \eqref{combipoint}, and
where the latter partition function $Z_{n}^{6V-U}[q^2,-1;q^{-2}]$ corresponds to taking uniform values $z_i=q^2$ and 
$w_i=w=q^{-8}=-1$ for all $i$, and $p=q^{-2}$
in the U-turn boundary partition function $Z_{n}^{6V-U}[\bz,\bw;p]$.

\begin{remark}\label{notcombirem}
We note that neither $Z_{n}^{6V-U}[q^2,-1;q^{-2}]$ nor $Z_{n}^{6V}[q^2,-1]$ are ``combinatorial" in the sense that they count objects. Indeed, they are weighted countings of configurations of the corresponding models, including irrational weights (such as $\sqrt{2}$).
For instance, we may choose to normalize the weights of the $6V$ model by an overall factor $\rho_o=q \sqrt{2}$ 
on odd rows such that $\rho_o(A,B,C)(qz,q^{-1}w)=(1,\sqrt{2},1)$ and $\rho_e=q^{-1} \sqrt{2}$ on even rows, such that $\rho_e(A,B,C)(qz^{-1},q^{-1}w)=(\sqrt{2},1,1)$ for $z=q^2$ and $w=-1$. In that case, the $6V$ partition function ${\bar Z}_n^{6V}$ with these new weights reads ${\bar Z}_n^{6V}=\rho_o^{n^2}\,\rho_e^{n(n-1)}\, Z_n^{6V}[q^2,-1]$. Yet, using the above relations at the combinatorial point, we find that 
${\bar Z}_n^{6V} =2^{n(n-1)/4}\, Z_n^{20V}$. So we should think of the $6V$ partition functions at hand as combinatorial tools rather than counting combinatorial objects.
\end{remark}

 \subsection{Inhomogeneous determinant formula}
 \label{kupsec}
 The result for the partition function of the U-turn boundary 6V model \cite{tsu,kuperberg2002symmetry} reads as follows for our choice of weights.
 For short, we denote by $(a,b,c)(z,w):=(a_3,b_3,c_3)(z,w)=(A,B,C)(q z ,q^{-1} w)$.
Let $M_U(n)$ be the $n\times n$ matrix with entries
\begin{equation}
[M_U(n)]_{i,j}=\frac{1}{a(z_i,w_j)b(z_i,w_j)}-\frac{1}{a(1 ,z_iw_j)b(1,z_iw_j)}
\end{equation}
Then the inhomogeneous partition function is given by:
\begin{eqnarray}\label{Udelta}
&&Z_{n}^{6V-U}[z_1,z_2,...,z_{n};w_1,w_2,...,w_n;p] \nonumber \\
&&\qquad=\det\big(M_U(n)\big)\times{
\scriptstyle 
\frac{\left\{\prod_{i=1}^n (p-p^{-1}w_i)a(q z_i,q^{-1}z_i^{-1})c(z_i,w_i) \right\}\left\{\prod_{i,j=1}^n a( z_i, w_j)b( z_i, w_j)
a(1 , z_iw_j)b(1,z_iw_j)\right\}}{\left\{\prod_{i=1}^n z_i^{n-1}\right\}\left\{\prod_{1\leq i<j \leq n}(z_i-z_j)(w_j-w_i)\right\}\left\{\prod_{1\leq i\leq j\leq n}(1-z_iz_j)(1-w_iw_j)\right\}}}
\end{eqnarray}

The aim of the following sections is to extract a homogeneous limit from this, with relevant values of spectral parameters, to finally get a compact formula for the partition function of the $20V$ model.

 \subsection{Homogeneous limit}
 \label{homsec}
 
 Let us denote by $\Delta_n[z_1,...,z_n;w_1,...,w_n]$ the quantity:
 \begin{equation}\label{defdelt}
 \Delta_n[z_1,...,z_n;w_1,...,w_n]:=\frac{ \det\big(M_U(n)\big)}{\left\{\prod_{1\leq i<j \leq n}(z_i-z_j)(w_j-w_i)\right\}\left\{\prod_{1\leq i\leq j\leq n}(1-z_iz_j)(1-w_iw_j)\right\}}
 \end{equation}
In view of the determinant formula \eqref{Udelta}, 
to deduce the homogeneous limit $Z_{2n\times n}^{6V-U}[z,w;p]$ \eqref{homog}, we must now compute the limit
$$\Delta_n[z,w]:=\lim_{z_i\to z \atop w_i\to w} \Delta_n[z_1,...,z_n;w_1,...,w_n] $$

In fact, we are only interested in the special combinatorial point values $z=q^2,w=q^8=-1$ of the homogeneous parameters.
In this section, we start by computing $\Delta_n[z,-1]$.

\noindent{\bf Notation:} Throughout the paper, for any function $f(u,v)$ with a power series expansion around $(0,0)$,
we denote by $f(u,v)\vert_{u^iv^j}$ the coefficient of $u^iv^j$ in the series expansion of $f$, in other words 
$f(u,v)\vert_{u^iv^j}:=\frac{\partial_u^i}{i!}\frac{\partial_v^j}{j!}f(0,0)$.

\begin{thm}\label{thmdelta}
We have:
\begin{equation}\label{formudel}
\Delta_n[z,-1]=\frac{(-1)^{n(n-1)/2}}{2^n (1-z^2)^{n(n+1)/2}}\det_{0\leq i,j\leq n-1}\left(
f_U(z;u,v)\Big\vert_{u^iv^{2j+1}} \right)  ,\end{equation}
where $f_U(z;u,v)$ denotes the function:
\begin{eqnarray}\label{fu}
f_U(z;u,v)&:=& \frac{1}{a(z+u,v-1)b(z+u,v-1)}-\frac{1}{a(1 ,(z+u)(v-1))b(1,(z+u)(v-1))} \nonumber \\
&=& {\scriptstyle \frac{1}{(q^{-2}-q^2)(z+u)}\left\{\frac{q^{-1}}{a(z+u,v-1)}-\frac{q}{b(z+u,v-1)}-\frac{q^{-1}}{a(\frac{1}{z+u} ,v-1)}+\frac{q}{b(\frac{1}{z+u},v-1)}\right\}}  
\end{eqnarray}
\end{thm}
\begin{proof}
The theorem is proved in two steps. We first keep $z_i$ generic, and compute for $w=-1$:
$$\Delta_n[z_1,...,z_n;-1]:=\lim_{w_i\to -1\atop i=1,2,...,n} \frac{ \det\big(M_U(n)\big)}{\left\{\prod_{1\leq i<j \leq n}(w_j-w_i)\right\}\left\{\prod_{1\leq i\leq j\leq n}(1-w_iw_j)\right\}}$$
Setting $f[z;v]:=f_U(z;0,v)$, and setting $w_i=-1+v_i$, we may rewrite:
$$\Delta_n[z_1,...,z_n;-1]:=\lim_{v_i\to 0\atop i=1,2,...,n} \frac{ \det_{1\leq i,j\leq n}\big(f[z_i;v_j]\big)}{\left\{\prod_{1\leq i<j \leq n}(v_j-v_i)\right\}\left\{\prod_{1\leq i\leq j\leq n}(v_i+v_j)\right\}}$$
where in the last factors $1-w_iw_j=v_i+v_j-v_iv_j$ we have dropped the higher order term $v_iv_j$ as the $v$'s tend to $0$.

For each $i$, we may expand $f[z_i;v_j]$ in power series of $v_j$:
$f[z_i;v_j]=\sum_{m\geq 0} f_m[z_i] v_j^m $, where the coefficients read:
\begin{equation}
\label{columns}
f_m[z_i]={\scriptstyle \frac{1}{(q^{-2}-q^2)z_i}\left\{ \left(\frac{1}{q^2 z_i+1}\right)^{m+1}+\left(\frac{q^2 z_i}{q^2z_i+1}\right)^{m+1}
-\left(\frac{1}{q^{-2} z_i+1}\right)^{m+1}-\left(\frac{q^{-2} z_i}{q^{-2}z_i+1}\right)^{m+1}\right\}}
\end{equation}
Let $\bF[v]$ (resp. $\bF_m$) denote the column vector with entries $f[z_i;v]$ (resp. $f_m[z_i]$) $i=1,2,...,n$. We have
$\bF[v]=\sum_{m\geq 0}v^m\,\bF_m $.
Moreover, the columns $\bF_m$ are not linearly independent. More precisely we have the following.

\begin{lemma}\label{lemcomb}
Let $c_n$, $n\in \Z_+$ be the coefficients in the series expansion:
\begin{equation}\label{bernou}
\frac{2}{1+e^{-x}}=\sum_{n\geq 0} c_n\, \frac{x^n}{n!}
\end{equation}
namely $c_0=1$, $c_1=\frac{1}{2}$, $c_2=0$, $c_3=-\frac{1}{4}$, $c_4=0$, $c_5=\frac{1}{2}$, $c_6=0$, $c_7=\frac{17}{8}$, ...
Then for all $n\geq 0$ we have:
\begin{equation}\label{lincomb}\bF_{2n}=\sum_{i=0}^{n-1} c_{2i+1}{2n+1\choose 2i+1} \bF_{2n-2i-1}
\end{equation}
\end{lemma}
\begin{proof}
For any function $f(\al)$ (say polynomial), introduce the shift operator $S: f\mapsto Sf$, such that $(Sf)(\al)=f(\al-1)$.
Using the Taylor expansion around $\al=0$, we may express $S$ as $\exp(-\frac{d}{d\al})$. In particular,
the identity $(1+\exp(-\frac{d}{d\al}))f (\al)= f(\al)+f(\al-1)$ is easily inverted using the series expansion \eqref{bernou} in which $x$ is substituted with the differential operator $\frac{d}{d\al}$:
$$f (\al) =\frac{1}{2} \sum_{n\geq 0} c_n \frac{1}{n!}\frac{d^n}{d\al ^n} \big(f(\al)+f(\al-1)\big) $$
Let us apply this to the polynomial $f(\al)=\al^{2n+1}$ for $n\geq 0$. We get:
$$\al^{2n+1} =\frac{1}{2} \left\{ \al^{2n+1}+(\al-1)^{2n+1} +\sum_{i=0}^{n} c_{2i+1} {2n+1\choose 2i+1} (\al^{2n-2i}+(\al-1)^{2n-2i})
\right\}$$
which implies:
$$\al^{2n+1}+(1-\al)^{2n+1}=\sum_{i= 0}^{n} c_{2i+1} {2n+1\choose 2i+1} (\al^{2n-2i}+(1-\al)^{2n-2i})\ .$$
Up to the overall prefactor $1/((q^{-2}-q^2)z_i)$, the $i$-th row of the desired identity \eqref{lincomb} follows immediately from this, as it takes precisely the form:
\begin{eqnarray*}&&\al^{2n+1}+(1-\al)^{2n+1}-(\beta^{2n+1}+(1-\beta)^{2n+1})\\
&&\quad =\sum_{i=0}^{n-1} c_{2i+1}{2n+1\choose 2i+1} ((\al^{2n-2i}+(\al-1)^{2n-2i})-(\beta^{2n-2i}+(\beta-1)^{2n-2i}))
\end{eqnarray*}
where  $\al=\frac{1}{q^2z_i+1}$ and $\beta=\frac{1}{q^{-2}z_i+1}$.
\end{proof}

Note that as $\bF_0=0$ we have the series expansion $\bF[v]=\sum_{n\geq 1} \bF_n v^n=:v\bG[v]$. 
The matrix with entries $f[z_i;v_j]$ is made of successive columns $v_1\bG[v_1],...,v_n\bG[v_n]$,
therefore the determinant of this matrix is $v_1v_2 \cdots v_n \det(\bG[v_1]\bG[v_2]\cdots \bG[v_n])$.
Consequently we rewrite:
$$\frac{ \det_{1\leq i,j\leq n}\big(f[z_i;v_j]\big)}{\left\{\prod_{1\leq i<j \leq n}(v_j-v_i)\right\}\left\{\prod_{1\leq i\leq j\leq n}(v_i+v_j)\right\}}
=\frac{1}{2^n} \frac{\det(\bG[v_1]\bG[v_2]\cdots \bG[v_n])}{\prod_{1\leq i<j \leq n}(v_j^2-v_i^2)}=:D_n[\bz;\bv]$$
Let us now successively take $v_1\to 0$, $v_2\to 0$, ... $v_n\to 0$. The first limit $v_1\to 0$ amounts to taking $\bG[0]=\bF_1$ as first column in the determinant and set $v_1=0$ in the denominator. To perform the second, we must expand $\bG[v_2]=\bF_1+\bF_2 v_2 + \bF_3v_2^2+O(v_2^3)$ up to order $2$, as both constant and linear terms are columns proportional to $\bF_1$ by Lemma \ref{lemcomb}. When we reach $v_k\to 0$, we must expand $\bG[v_k]$ up to order $2k$, as all previous terms are linear combinations of the previous columns in the determinant, as a consequence of  Lemma \ref{lemcomb}.
This leads  successively  to:
\begin{eqnarray*} \lim_{v_1\to 0}D_n[\bz;\bv]&=&\frac{1}{2^n}\frac{\det(\bF_1\bG[v_2]\cdots \bG[v_n])}{\prod_{i=2}^n v_i^2 \prod_{2\leq i<j \leq n}(v_j^2-v_i^2)}\\
\lim_{v_1,v_2\to 0}D_n[\bz;\bv]&=&\frac{1}{2^n}\frac{\det(\bF_1\bF_3\bG[v_3]\cdots \bG[v_n])}{\prod_{i=3}^n v_i^4 \prod_{3\leq i<j \leq n}(v_j^2-v_i^2)}\\
\quad \vdots \ && \quad \\
\lim_{v_1,v_2,...,v_n\to 0}D_n[\bz;\bv]&=&\frac{1}{2^n}\det(\bF_1\bF_3\cdots \bF_{2n-1})=\Delta_n[z_1,...,z_n;-1]
\end{eqnarray*}

We are now left with the second step, namely sending all $z_i$'s to $z$:
\begin{eqnarray*}\Delta_n[z,-1]&=&\lim_{z_i\to z\atop i=1,2,...,n} \frac{\Delta_n[z_1,...,z_n;-1]}{\prod_{1\leq i<j \leq n}(z_i-z_j)\prod_{1\leq i\leq j\leq n}(1-z_iz_j)}\\
&=&\lim_{u_i\to 0\atop i=1,2,...,n} \frac{1}{2^n(1-z^2)^{n(n+1)/2}} \frac{\det(\bF_1\bF_3\cdots \bF_{2n-1})}{\prod_{1\leq i<j \leq n}(u_i-u_j)}
\end{eqnarray*}
where we have set $z_i=z+u_i$, and
the column vector $\bF_{2j+1}$ has entries $f_{2j+1}[z+u_i]$, $i=1,2,...,n$. Let us represent the matrix 
$(\bF_1 \bF_3 \cdots \bF_{2n-1})$ by its rows, in the form $(\bR[u_1] \bR[u_2] \cdots \bR[u_n])^t$. Here, the column vector $\bR[u_i]$ is the transpose of the $i$-th row vector of the original matrix, and $\bR[u]$ has entries $f_{2j+1}[z+u]$, $j=0,1,...,n-1$.
As before, we write the series expansion $\bR[u]=\sum_{n\geq 0} \bR_n u^n$, and perform the successive limits $u_1\to 0$, $u_2\to 0$, ... $u_n\to 0$ in the quantity:
$$D_n[\bu]:=\frac{1}{2^n(1-z^2)^{n(n+1)/2}}\frac{\det( \bR[u_1] \bR[u_2] \cdots \bR[u_n])^t }{\prod_{1\leq i<j\leq n} (u_i-u_j)}$$
The vectors $\bR_i$ turn out to be linearly independent, so that each vector $\bR[u_j]$ only has to be expanded up to order $j-1$,
and we end up with the formulas: 
\begin{eqnarray*}
\lim_{u_1\to 0}D_n[\bu]&=& \frac{(-1)^{n-1}}{2^n(1-z^2)^{n(n+1)/2}}\frac{\det( \bR_0 \bR[u_2] \cdots \bR[u_n])^t }{\prod_{i=2}^n u_i \prod_{2\leq i<j\leq n} (u_i-u_j)}\\
\lim_{u_1,u_2\to 0}D_n[\bu]&=& \frac{(-1)^{n-1+n-2}}{2^n(1-z^2)^{n(n+1)/2}}\frac{\det( \bR_0 \bR_1 \bR[u_3]\cdots \bR[u_n])^t }{\prod_{i=3}^n u_i^2 \prod_{3\leq i<j\leq n} (u_i-u_j)}\\
\quad \vdots && \quad \\
\lim_{u_1,u_2,...,u_n\to 0}D_n[\bu]&=&\frac{(-1)^{n(n-1)/2}}{2^n(1-z^2)^{n(n+1)/2}}\det( \bR_0 \bR_1 \cdots \bR_n)^t=\Delta_n[z;-1]
\end{eqnarray*}
The identity \eqref{formudel} follows by identifying the entries $i+1,j+1$ ($i,j=0,1,...,n-1$) of the matrix $( \bR_0 \bR_1 \cdots \bR_n)^t$ with the coefficient of $u^i v^{2j+1}$ of the series expansion of $f_U(u,v)$ around $(0,0)$. Indeed, $(\bF_{2j+1})_{i+1}$ corresponds to the coefficient of $v^{2j+1}$ in the series expansion of $\bF[v]$ around $v=0$,
and $(\bR_i)_{j+1}$ to the coefficient of $u^i$ in the expansion of $(\bF_{2j+1})_{i+1}$ around $u=0$, after setting $z_i=z+u$.
This completes the proof of the theorem.
\end{proof}

To further the computation of $\Delta_n[z;-1]$, we need the following.
\begin{lemma}\label{lemdeltatwo}
Let ${\tilde f}_U(u,v):=\gamma(u)\,\delta(v)  \,f_U\Big(a u\al(u),b v\beta(v)\Big)$, with $a,b\in \C^*$, and where $\al(x),\beta(x),\gamma(x),\delta(x)$ are power series that converge for small enough $x$, and take the value $1$ at $x=0$.
Then we have:
$$\Delta_n[z,-1]=\frac{(-1)^{n(n-1)/2}}{2^n (1-z^2)^{n(n+1)/2}a^{n(n-1)/2} \,b^{n^2} }\det_{0\leq i,j\leq n-1}\left(
{\tilde f}_U(u,v)\Big\vert_{u^iv^{2j+1}} \right) $$
\end{lemma}
\begin{proof}
Repeating the first step of the proof of Theorem \ref{thmdelta} with ${\tilde f}_U$ substituted for $f_U$, we see that when we send $v_j\to 0$, the column 
$\bF_{2j+1}$ gets multiplied by $b^{2j+1}$, but otherwise the substitution $v_j\to b v_j \beta(v_j)$, 
as well as the overall factor $\delta(v_j)$ only affect higher order terms of the expansion in $v_j$, 
hence the limiting process yields the same result as before, up to an overall factor $b^{1+3+\cdots +(2n-1)}=b^{n^2}$. 
The same happens for the second step, where each vector $\bR_i$ gets multiplied by $a^i$ leading to an overall factor 
$a^{n(n-1)/2}$, while the substitution $u_i\to a u_i\al(u_i)$, as well as the overall factor $\gamma(u_i)$
again only affect higher order terms in the expansion in $u_i$. The Lemma follows.
\end{proof}

We now compute the quantity $\Delta_n[q^2,-1]$, corresponding to $z=q^2$.
\begin{thm}\label{deltathm}
We have 
\begin{eqnarray}
\Delta_n[q^2;-1]&=&
\frac{
q^{6n(n+1)} }{2^{\frac{n(7n+5)}{4}}}\det_{0\leq i,j\leq n-1}\left(
g_U(u,v)\Big\vert_{u^iv^j} \right) \label{deltaunif}\\
g_U(u,v)&=& \frac{(1+u^2)(1+2u-u^2)}{(1-u^2v)\Big((1-u)^2-v(1+u)^2\Big)}
\label{gu}
\end{eqnarray}
\end{thm}
\begin{proof}
We apply Lemma \ref{lemdeltatwo} with $z=q^2$, $a=2q^6$, $b=2 q^{-4}$, $\al(u)=\frac{1}{1-q^4 u}$, $\beta(v)=\frac{1}{1+q^{-4}v}$,
$\gamma(u)=\frac{1+q^4 u}{1-q^4 u}$, $\delta(v)=\beta(v)^2$, with the result:
\begin{eqnarray*}{\tilde f}_U(u,v)&=&\frac{1+q^4 u}{(1-q^4 u)(1+q^{-4}v)^2} 
f_U\left(\frac{2q^6  u}{1-q^4 u}, \frac{2 q^{-4} v}{1+q^{-4}v}\right)  \\
&=&v\, \frac{q^4}{1-q^4}\, \frac{(1+u^2)(1+2u-u^2)}{(1-u^2v^2)\Big((1-u)^2-v^2(1+u)^2\Big)}\\
&=&-\frac{q^{-2}}{\sqrt{2}} \, v\, g_U(u,v^2)
\end{eqnarray*}
%
Noting that the coefficient of $v^{2j+1}$ in ${\tilde f}_U(u,v)$ corresponds to the coefficient of $v^j$ in $g_U(u,v)$,
the theorem follows from assembling all the factors and using the fact that $q=e^{i\pi/8}$.
\end{proof}

\subsection{The 20V model partition function}
\label{res20vsec}
We now complete the calculation of the partition function for the $20V$ model on the quadrangle $\cQ_n$.

\begin{thm}\label{thm20V}
The partition function of the 20V model with uniform weights on the quadrangle ${\mathcal Q}_n$ reads:
\begin{eqnarray}
Z_n^{20V}&=&\det_{0\leq i,j\leq n-1}\left(
g_{20V}(u,v)\Big\vert_{u^iv^j} \right)\label{z20v} \\
g_{20V}(u,v)&=&\frac{(1+u^2)(1+2u-u^2)}{(1-u^2v)\Big((1-u)^2-v(1+u)^2\Big)}=\frac{(\frac{1+u}{1-u})^2}{1-v(\frac{1+u}{1-u})^2} -\frac{u^2}{1-v u^2}\label{f20v}
\end{eqnarray}
\end{thm}
\begin{proof}
We use the result of Theorem \ref{deltathm} to compute the homogeneous limit of the determinant formula \eqref{Udelta} 
for $Z^{6V-U}_{n}[\bz,\bw;p]$, with $z_i\to q^2$, $w_i\to -1$ and $p=q^{-2}$.
First we assemble all the
remaining factors in $Z^{6V-U}_{n}[\bz,\bw;p]$ to find that:
\begin{eqnarray*}Z^{6V-U}_{n}[q^2,-1]&=&
\frac{1}{z^{n(n-1)}}\left((q^2-q^{-2})\sqrt{z w}\right (p^{-2}-p^{2})(p-p^{-1}w))^n\\
&\times& \left((q z-q^{-1}w)(q^{-1}z -q w)(q-q^{-1}z w)(q^{-1}-q z w)\right)^{n^2} \, \Delta_n[z,w]\vert_{z\to q^2,w\to -1}\\
&=& q^{-2n(n-1)}(q^2-q^{-2})^n(q^4-q^{-4})^n (q^2+q^{-2})^n  (8 q^4)^{n^2} \Delta_n[q^2;-1]\\
&=& q^{n(2n-1)} 2^{n(3n+2)} (-1)^n   \Delta_n[q^2;-1]
\end{eqnarray*}
Substituting the result of Theorem \ref{deltathm} for $\Delta_n[q^2;-1]$, we finally get:
\begin{equation*}
Z^{6V-U}_{n}= q^{5n} 2^{\frac{n(5n+3)}{4}} \det_{0\leq i,j\leq n-1}\left(
g_U(u,v)\Big\vert_{u^iv^j} \right)
\end{equation*}
The theorem follows from combining this with \eqref{20vone}, and defining $g_{20V}(u,v):=g_U(u,v)$.
\end{proof}

\begin{remark} \label{rem20V}
We may interpret this result as follows. Let $P$ be the infinite matrix with entries $P_{i,j}$, $i,j\in \Z_+$ generated by the series expansion
$g_{20V}(u,v)=\sum_{i,j\geq 0} u^i v^j P_{i,j}$ of \eqref{f20v}. 
Then $Z^{20V}_n=\det(P_n)$ where $P_n$ is the truncation of the matrix $P$
to it $n$ first rows and columns.
\end{remark}

\section{Domino tilings of the Aztec triangle}
\label{dtsec}

\subsection{Path model}
\label{pathsec}

\begin{figure}
\begin{center}
\includegraphics[width=13cm]{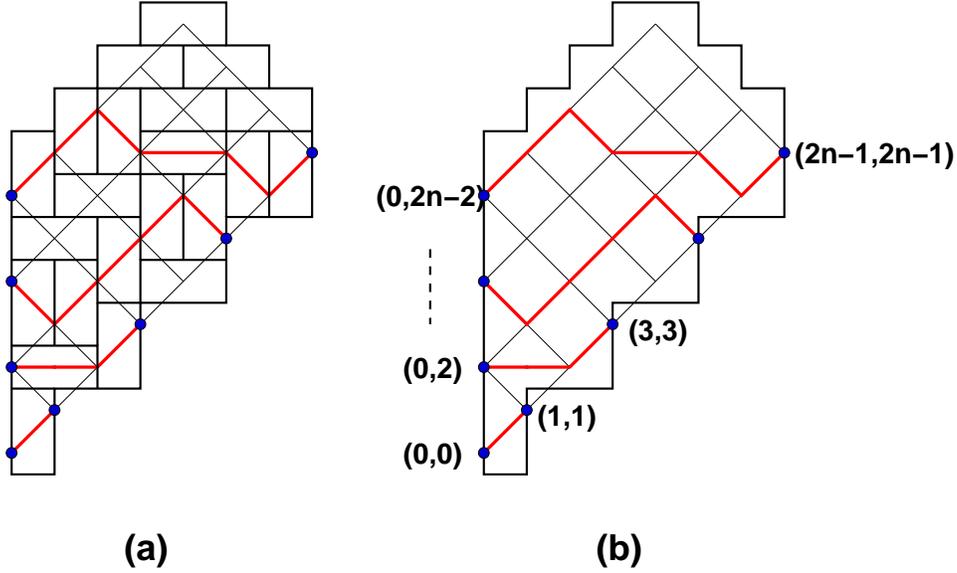}
\end{center}
\caption{\small The Aztec triangle $\cT_n$ and a typical domino tiling configuration (a), together with its bijectively associated family of non-intersecting Schr\"oder paths. In (b) we show the starting points of the $n$ paths: $(0,2i)$, $i=0,1,...,n-1$ and their endpoints $(2j+1,2j+1)$, $j=0,1,..,n-1$.
}
\label{fig:pachtermax}

\end{figure}

The Aztec triangle $\cT_n$ is defined as the planar domain of Figs.~\ref{fig:pachtermax} (a,b),
where the length of the vertical West boundary is $2n$. We study and count the 
tiling configurations of this domain by means of $2\times 1$ dominos. The latter are in bijection with configurations
of $n$ non-intersecting lattice paths from the West boundary (starting points $(0,2i)$, $i=1,2,...,n-1$) to the staircase
East boundary (endpoints $(2j+1,2j+1)$, $j=0,1,...,n-1$). 
These paths are Schr\"oder paths, with three types of allowed steps on $\Z^2$: up $(1,1)$, down $(1,-1)$ and horizontal $(2,0)$.

Let $Z_{n}^{DT}$ denote the total number of domino tiling configurations of the Aztec triangle $\cT_n$.
The latter is easily computed by direct application of the Lindstr\"om-Gessel-Viennot determinant formula \cite{LGV1,GV}.
Denoting by $M_{i,j}$ the number of configurations of a single path from the starting point $(0,2i)$ to the endpoint
$(2j+1,2j+1)$, we have:
\begin{lemma}\label{thmzat}
The partition function for domino tilings of the Aztec triangle $\cT_n$ reads:
$$Z_{n}^{DT}=\det_{0\leq i,j \leq n-1} \left(M_{i,j} \right) $$
\end{lemma} 

\begin{remark}\label{remAT}
Note that $M_{i,j}$ is independent of $n$, as the shape of the Aztec triangle domain imposes no further constraint 
on the paths apart from their fixed starting and endpoints. The Lemma \ref{thmzat} states that $Z_{n}^{DT}=\det(M_n)$, where 
the $n\times n$ matrix $M_n$ may be viewed as the {\it finite} truncation 
to the first $n$ rows and columns of the infinite matrix  $M$ with entries $M_{i,j}$, $i,j\in \Z_+$.
\end{remark}

\subsection{Generating functions and determinant formula}
\label{resdtsec}

The generating function $f_M(u,v):= \sum_{i,j=0}^\infty v^i u^j \, M_{i,j}$ is given by the following.

\begin{thm}\label{genDT}
We have
$$ f_M(u,v)= \sum_{i,j=0}^\infty v^i u^j \, M_{i,j} = \frac{1+u}{1-v-4 u v-u^2 v+u^2 v^2} $$
\end{thm}
\begin{proof}
We start from the generating series:
\begin{equation}\label{defsig} \sigma(u,v)=\frac{1}{1-u-v-uv}
\end{equation}
We may interpret $\sigma(u,v)$ as the partition function for arbitrary finite length Schr\"oder paths from the origin, with weight $u$ per up step, $v$ per down step, and $uv$ per horizontal step. The position of the endpoint of the path
reads $(U+D+2H,U-D)$ where $U,D,H$ denote respectively the total numbers of up, down, horizontal steps in the path.
Such a path receives a weight $u^U v^D (u v)^H=u^{U+H} v^{D+H}$. The number $\Sigma_{i,j}$ of paths
from the origin with fixed endpoint $(i,j)$, with $i-j$ even, therefore reads:
$\Sigma_{i,j}=\sigma(u,v)\vert_{u^{\frac{i+j}{2}} v^{\frac{i-j}{2}}}$. Equivalently:
\begin{equation}\label{sigentry}
\Sigma_{i,j}=\oint \frac{dx}{2i\pi x} \frac{dy}{2i\pi y}\, \frac{1}{x^{\frac{i+j}{2}} y^{\frac{i-j}{2}}} \,\frac{1}{1-x-y-x y}
\end{equation}
where the contour integral picks the residues at $0$. By a trivial translation of the starting point, we may now interpret $M_{i,j}$ as the total number of paths from the origin to $(2j+1,2j+1)-(0,2i)=(2j+1,2j+1-2i)$, hence
$M_{i,j}=\Sigma_{2j+1,2j+1-2i}$, or:
$$M_{i,j}=\oint \frac{dx}{2i\pi x} \frac{dy}{2i\pi y}\, \frac{1}{x^{2j+1-i} y^{i}} \, \frac{1}{1-x-y-x y}$$
Changing integration variables to $u=y/x, v=x$, we get
$$M_{i,j}=\oint \frac{du}{2i\pi u} \frac{dv}{2i\pi v}\, \frac{1}{v^{2j+1} u^{i}} \, \frac{1}{1-v-u v-u v^2}=\frac{1}{1-v-u v-u v^2}\Big\vert_{u^iv^{2j+1}}$$
Finally, we remark that for any series $f(v)$, and any $j\in \Z_+$, we have:
\begin{equation}\label{trick}f(v)\vert_{v^{2j+1}}=\frac{1}{2 \sqrt{v}}\left( f(\sqrt{v})-f(-\sqrt{v}) \right) \Big\vert_{v^j} 
\end{equation}
and the theorem follows from the identity
$$\frac{1}{2\sqrt{v}}\left(\frac{1}{1-(1+u)\sqrt{v}-u v}-\frac{1}{1+(1+u)\sqrt{v}-u v}\right)=\frac{1+u}{1-v-4 u v-u^2 v+u^2 v^2}
$$
\end{proof}

This allows to translate the result of Lemma \ref{thmzat} into the following determinant formula.

\begin{thm}\label{thmAT}
The partition function for domino tilings of the Aztec triangle $\cT_n$ reads:
\begin{eqnarray}
Z^{DT}_n&=&\det_{0\leq i,j \leq n-1}\left( 
f_{DT}(u,v)\Big\vert_{u^iv^j} \right)\label{zat} \\
f_{DT}(u,v)&=& \frac{1+u}{1-v -4 u v - u^2 v +u^2 v^2} \label{fat}
\end{eqnarray}
\end{thm}


\section{20V-Domino Tiling correspondence}
\label{correspsec}

\subsection{Determinant relations}\label{detrelasec}

We wish to compare the results of Theorems \ref{thm20V} and \ref{thmAT}.
Using Remarks \ref{rem20V} and \ref{remAT}, we see that both enumeration results are expressed as the determinant of a finite
truncation of an infinite matrix. We have the following \cite{DFLAP}.

\begin{lemma}\label{detlem}
Let $L$ (resp. $U$) be arbitrary infinite lower (resp. upper) triangular matrices, with diagonal entries $1$, and $A$ any infinite matrix. 
Then we have the following properties:
\begin{enumerate}
\item For all $n\geq 1$, the truncation $(LAU)_n$ of the infinite matrix $LAU$ to its first $n$ rows and columns is equal to $L_n A_n U_n$, the product of the corresponding truncations of $L,A,U$ respectively.
\item $\det(A_n)=\det\big((LAU)_n\big)$
\end{enumerate}
\end{lemma}
\begin{proof}
We have for all $i,j \in [0,n-1]$:
$((LAU)_n)_{i,j}=\sum_{\ell,m\geq 0} L_{i,\ell} A_{\ell,m}U_{m,j}$. But $L_{i,\ell}=0$ for $\ell>i$ hence in particular for all $\ell\geq n$, as $i\leq n-1$. Similarly, $U_{m,j}=0$ for $m\geq n$, and therefore we may rewrite
$((LAU)_n)_{i,j}=\sum_{\ell,m=0}^{n-1} L_{i,\ell} A_{\ell,m}U_{m,j}=((L_nA_nU_n)_{i,j}$, and (1) follows. The statement (2) follows from
$\det(L_nA_nU_n)=\det(A_n)$ as $\det(L_n)=\det(U_n)=1$.
\end{proof}

\begin{lemma}\label{prodlem}
Let $A,B$ be two infinite matrices with entries generated by the series $f_A(u,v)=\sum_{i,j\geq 0} u^i v^j A_{i,j}$
and $f_B(u,v)=\sum_{i,j\geq 0} u^i v^j B_{i,j}$. Then the product $AB$ has for generating series the convolution of $f_A$ and $f_B$,
namely:
$$f_{AB}(u,v)=f_A*f_B(u,v)=\oint \frac{dt}{2i\pi t} f_A(u,t^{-1})f_B(t,v) $$
where the contour integral picks up the residue at $t=0$.
\end{lemma}
\begin{proof} By direct calculation.
\end{proof}

\subsection{Proof of the equivalence theorem}
\label{equivsec}

\begin{thm}\label{equivalencethm}
The partition function of the 20V DWBC3 model on the quadrangle ${\mathcal Q}_n$ and that of the domino tilings of the Aztec triangle 
$\cT_n$ coincide for all $n\geq 1$.
\end{thm}
\begin{proof}
By Remarks \ref{rem20V} and \ref{remAT}, we have $Z_n^{20V}=\det(P_n)$, $Z_n^{DT}=\det(M_n)$ in terms of truncations of infinite
square matrices $P,M$ respectively generated by $f_P=g_{20V}$ \eqref{f20v} and $f_M=f_{DT}$ \eqref{fat}.
Let $L$ be the infinite matrix with generating series 
\begin{equation}\label{defL}f_L(u,v):=\frac{1+2u-u^2}{1-u(1+v+uv)}\ .\end{equation}
Then we have the following properties: (1) $L$ is lower triangular, with diagonal elements $1$ and (2) $LM=P$.
The property (1) follows from noticing that we have a series identity of the form $f_L(u,v)=h(u,uv)$ where
$$h(u,v)=\frac{1+2u-u^2}{1-u-v-uv}\ , $$
which shows that non-zero coefficients of $u^i v^j$ always have $j\geq i$. Moreover the terms with $i=j$ correspond to 
the series expansion of $h(0,uv)=1/(1-uv)$ and are all equal to $1$. Property (2) follows from Lemma \ref{prodlem}:
\begin{eqnarray*}f_{LM}(u,v)&=&\oint \frac{dt}{2i\pi } \,\frac{1+2u-u^2}{t(1-u)-u(1+u)} f_M(t,v) \\ 
&=&\frac{1+2u-u^2}{1-u} f_M\left(u \frac{1+u}{1-u},v\right) 
= f_P(u,v) 
\end{eqnarray*}
easily checked by substituting the values $f_P=g_{20V}$ \eqref{f20v} and $f_M=f_{DT}$ \eqref{fat}.
Finally, by Lemma \ref{detlem} we have $\det(P_n)=\det((LM)_n)=\det(L_n M_n)=\det(M_n)$, and the theorem follows.
\end{proof}

\section{Refined partition functions of the 20V model}
\label{ref20Vsec}

\subsection{Setting and weights}

In this section, we consider the particular inhomogeneous partition function $Z_n^{20V}[w]$ for the 20V model on the quadrangle of Fig.~\ref{fig:domain}, corresponding to a choice of uniform spectral parameters
$z_i=q^2 z=q^4$, $i=1,2,...,n$, $t_i=z^{-1}=q^{-2}$, $i=1,2,...,2n-1$ and $w_i=q^8=-1$, $i=1,2,...,n-1$, while $w_n=w$ is kept arbitrary.
The corresponding weights \eqref{weights20V} are all equal to $1$, except those of the last column,
which read:
\begin{eqnarray}
 {\bar\omega}_0&=&\frac{(1-w)(q^2-q^{-2} w)}{2\sqrt{2}}\nonumber \\ 
 {\bar\omega}_1&=&\frac{(1-w)(q^{-2}-q^{2} w)}{2\sqrt{2}}\nonumber \\ 
{\bar\omega}_2={\bar\omega}_4&=&\frac{1-w}{2}\, \sqrt{-w}\nonumber \\ 
{\bar\omega}_3={\bar\omega}_5={\bar\omega}_6&=&\left(\frac{1-w}{2}\right)^2 \label{barog}
\end{eqnarray}
The quantity  $Z_n^{20V}[w]$ gives access to refined 20V model partition functions as explained in the next section.

\subsection{Refined partition functions}

\begin{figure}
\begin{center}
\includegraphics[width=8cm]{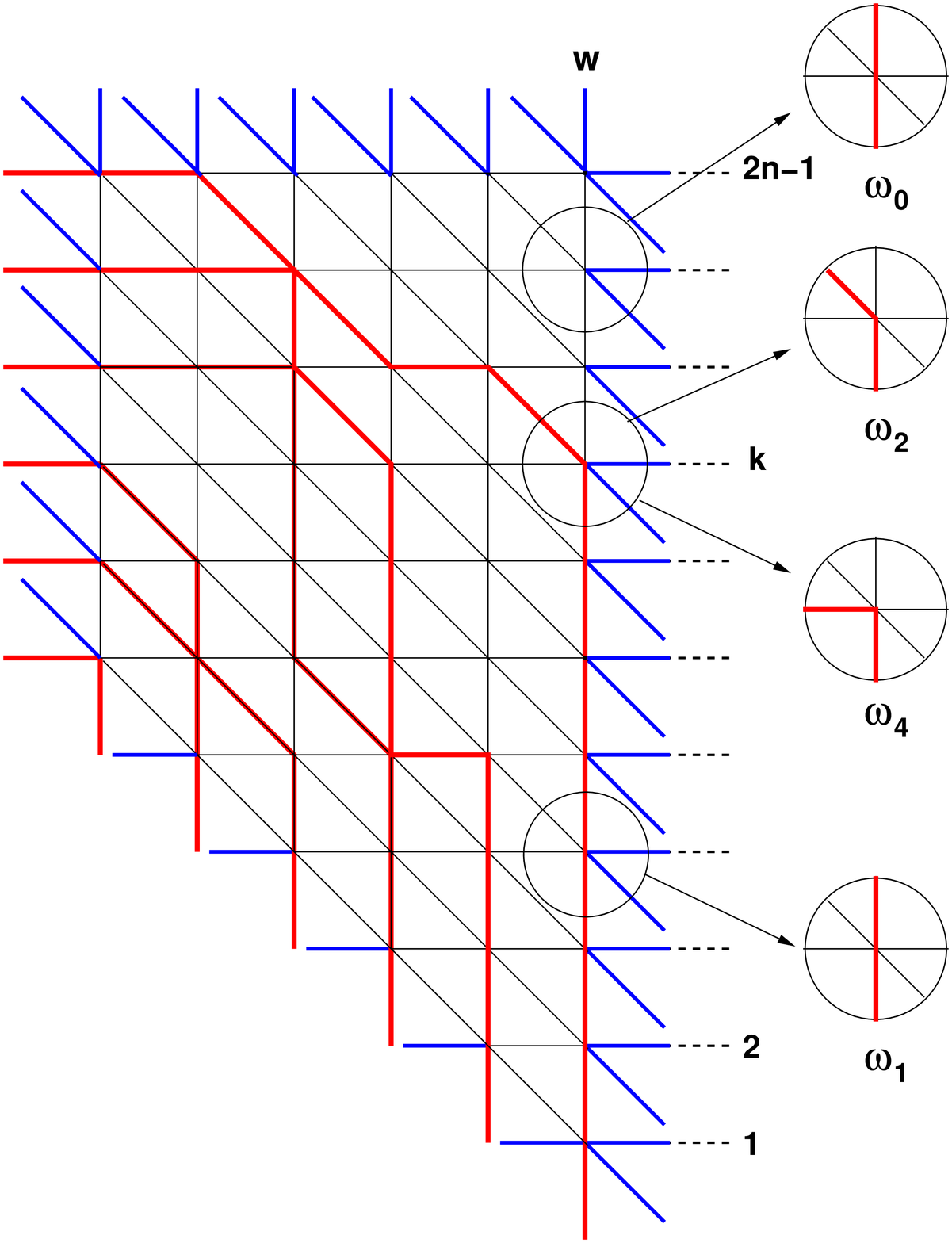}
\end{center}
\caption{\small A typical configuration of osculating Schr\"oder paths on the quadrangle $\cQ_n$. The vertex labeled $k\in [1,2n-1]$ corresponds to the first visit of any path to the last column. All vertices strictly above it are in the empty configuration (weight $\omega_0$), while those strictly below have a single vertical path going through (weight $\omega_1$). In general configurations, the vertex $k$ itself could be in either of the two indicated states (cameos on the right, with weights $\omega_2$ or $\omega_4$).
}
\label{fig:ref20v}

\end{figure}

Let us now relate the partition function $Z_n^{20V}[w]$ to refined configuration numbers of the 20V model on the quadrangle ${\cQ}_n$.
The osculating Schr\"oder path configurations of the 20V model on $\cQ_n$ may be classified according to the 
configuration of the vertices in the last column (see Fig.~\ref{fig:ref20v} for an illustration). Let $k$ denote the position (labeled $1,2,...,2n-1$ from the bottom) of the unique vertex of the last column where 
the top vertical edge is empty and the bottom vertical edge is occupied by a path. Alternatively, this $k$-th vertex corresponds to the point of entry of the topmost path into the rightmost vertical line, below which the path goes down vertically until its endpoint. Note that all vertices strictly above $k$
are in the empty configuration (weights $\omega_0$), and all those strictly below in the configuration with exactly one vertical line passing through (weights $\omega_1$). The vertex $k$ itself may be in either of two configurations, corresponding to whether the path enters the last column horizontally ($-$, weight $\omega_4$) or diagonally ($\backslash$, weight $\omega_2$).
Let $Z_{n,k}^{20V -}$ and $Z_{n,k}^{20V \backslash}$ denote the corresponding refined partition functions, with the obvious sum rule at $w=-1$:
$$\sum_{k=1}^{2n-1} Z_{n,k}^{20V -}+Z_{n,k}^{20V \backslash}=Z_{n}^{20V}$$
In the presence of the spectral parameter $w$ in the last column, this becomes:
\begin{equation}\label{sumrule}
\sum_{k=1}^{2n-1} \left({\bar\omega}_4 \,Z_{n,k}^{20V -}+{\bar \omega}_2\,Z_{n,k}^{20V \backslash}\right) {\bar \omega}_0^{2n-k-1} {\bar\omega}_1^{k-1}=Z_{n}^{20V}[w]
\end{equation}
with ${\bar\omega}_i$ as in \eqref{barog}.
Let us introduce a parameter: 
\begin{equation}
\label{deftau}
\tau:= \frac{{\bar\omega}_1}{{\bar\omega}_0}=\frac{q^{-2}-q^{2} w}{q^{2}-q^{-2} w}
\end{equation}
Further noting that ${\bar\omega}_4={\bar \omega}_2$, we get the following expression for the generating function $Z_n^{20V}(\tau)$
of refined 20V partition functions:
\begin{equation}\label{pfsumrule}
Z_n^{20V}(\tau):=\sum_{k=1}^{2n-1}Z_{n,k}^{20V}\,\tau^{k-1}=\frac{q^2-q^{-2}w}{{\bar \omega_0}^{2n-1}\sqrt{-2w}}\, Z_{n}^{20V}[w]
\end{equation}
where we used the notation $Z_{n,k}^{20V}=Z_{n,k}^{20V -}+Z_{n,k}^{20V \backslash}$ for the refined 20V partition function
corresponding to the leftmost path entering the rightmost vertical line at position $k$.
It is customary to also introduce the generating function $h_n^{20V}(\tau)$ for normalized ``one-point functions" 
$\frac{Z_{n,k}^{20V}}{Z_{n}^{20V}}$:
\begin{equation}\label{sumrule2}
h_n^{20V}(\tau):=\sum_{k=1}^{2n-1}\frac{Z_{n,k}^{20V}}{Z_{n}^{20V}}\,\tau^{k-1}=\frac{Z_{n}^{20V}(\tau)}{Z_{n}^{20V}}
=\frac{q^2-q^{-2}w}{{\bar \omega_0}^{2n-1}\sqrt{-2w}}\, \frac{Z_{n}^{20V}[w]}{Z_{n}^{20V}} 
\end{equation}

\subsection{Transformations to the 6V model}

Repeating the operations of Sections \ref{transfosec} and \ref{usec}, and recording all the weights of trivial vertices, we arrive at:
\begin{equation}\label{boil}
Z_n^{20V}[w]= {\scriptstyle q^{n^2}b_1(q^4,w)^{n-1} a_3(q^{-2},w)^n a_2(q^4,q^{-2})^{\frac{n(3n-1)}{2}}  b_1(q^4,-1)^{\frac{(n-1)(n-1)}{2}} \, a_3(q^{-2},-1)^{\frac{n(n-1)}{2}}}\, Z_{n}^{6V}[w]  
\end{equation}
where $Z_{n}^{6V}[w]$ denotes the partition function function of the 6V model with boundary conditions as in Fig.~\ref{fig:Uboundary} (c), with the vertex weights $(A,B,C)(qz_i,q^{-1} w_j)$
identical to those of the sublattices 2 and 3, with vertical spectral parameters $w_1=\cdots =w_{n-1}=-1$, $w_n=w$,
and horizontal spectral parameters alternating between the value $z=q^2$ on odd rows and $z^{-1}=q^{-2}$ 
on even rows.
Eq. \ref{boil} boils down to:
\begin{equation}\label{rela206}
\frac{Z_n^{20V}[w]}{Z_n^{20V}}= \left(\frac{1-w}{2}\right)^{2n-1} \, \frac{Z_{n}^{6V}[w]}{Z_{n}^{6V}} 
\end{equation}
trivially satisfied for $w=-1$.
Next,  \eqref{homog} for $p=q^{-2}$ implies:
\begin{equation}\label{rela6VU}
Z_{n}^{6V}[w]=  \frac{Z_{n}^{6V-U}[w]}{(p^2-p^{-2})^n(p+p^{-1})^{n-1}(p-p^{-1}w)}  
\end{equation}
and equivalently:
\begin{equation}
\frac{Z_n^{20V}[w]}{Z_n^{20V}}
=\frac{\sqrt{2}}{q^{-2}- q^{2} w}\,  \left(\frac{1-w}{2}\right)^{2n-1} \,\frac{Z_{n}^{6V-U}[w]}{Z_{n}^{6V-U}}
\end{equation}

\subsection{Refined partition functions of the  6V model}

\begin{figure}
\begin{center}
\includegraphics[width=9cm]{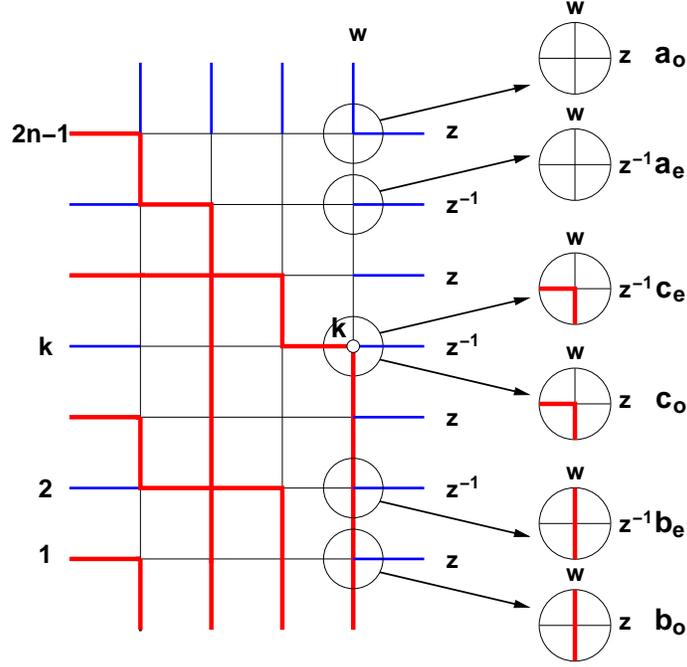}
\end{center}
\caption{\small A typical configuration of osculating paths for the $6V$ model on a rectangular grid of size $(2n-1)\times n$. The boundary conditions match those of  Fig.~\ref{fig:Uboundary} (c) upon global rotation by 180${}^\circ$. The horizontal spectral parameters alternates between the values $z=q^2$ on odd lines and $z^{-1}=q^{-2}$ on even lines. The vertical parameters are $w_i=q^8=-1$ except in the last column, where $w_n=w$.
The vertex labeled $k\in [1,2n-1]$ corresponds to the first visit of any path to the last column. All vertices strictly above it are in the empty configuration (weight $a_e,a_o$ on even/odd rows), while those strictly below have a single vertical path going through (weight $b_e,b_o$ on even/odd rows). In general configurations, the vertex $k$ itself could be in either of the two indicated states (cameos on the right, with weights $c_e$  if $k$ is even or $c_o$ if $k$ is odd).
}
\label{fig:Uref}
\end{figure}

We now turn to the partition function $Z_{n}^{6V}[w]$ of the 6V model. To easier compare it to that  of the 6V model 
with U-turn boundaries,
let us rotate the picture of Fig.~\ref{fig:Uboundary} (c) by 180${}^\circ$. The rows are now labeled $1,2,...,2n-1$ from bottom to top, and
the non-trivial spectral parameter $w$ is attached to the rightmost vertical line.
We consider the osculating path configurations of the model (see Fig.~\ref{fig:six} bottom), with $n$ paths starting on odd rows along the left boundary and ending 
on the S boundary (see Fig.~\ref{fig:Uref} for an illustration). We note that the topmost path must first visit the rightmost vertical line at a vertex corresponding to some row $k$,
and then continue along the vertical until its endpoint. The configurations along the rightmost vertical line have therefore b-type weights 
at all positions strictly under $k$, and a-type weights at all positions strictly above $k$, while the vertex at position $k$ is of type c.
When collecting the weights however, we must distinguish the cases according to the parity of $k$, as odd and even horizontal spectral parameters are different. Using the odd row weights: $(A,B,C)(q^3,q^{-1}w)$ and
even row weights $(A,B,C)(q^{-1},q^{-1}w)$, we get the following relative weights (ratios of the weights at generic $w$ to those at $w=-1$) in the last column:
$$ \begin{matrix} {\bar a}_o=\frac{q^2-q^{-2}w}{\sqrt{2}} & {\bar b}_o=\frac{1-w}{2} & {\bar c}_o=\sqrt{-w} \\
{\bar a}_e=\frac{1-w}{2} & {\bar b}_e=\frac{q^{-2}-q^{2}w}{\sqrt{2}} & {\bar c}_e=\sqrt{-w}\end{matrix}$$
where the subscript $e/o$ stands for even/odd row.

Let $Z_{n,k}^{6V}$ denote the refined partition function\footnote{As already mentioned in Remark \ref{notcombirem}, the 
quantities $Z_{n,k}^{6V}$ are not integers as the weighting of the $6V$ model is not uniform and involves $\sqrt{2}$ weights. Here and in the following we have dropped the arguments $[z=q^2,w=-1]$ indicating that we are at the combinatorial point, and we use the notation $[w]$ to indicate a non-trivial value $w$ for the $n$-th vertical spectral parameter, while all other spectral parameters have the values of the combinatorial point.} corresponding to paths first visiting the rightmost vertical at position $k$.
Apart from the trivial sum rule at $w=-1$:
$$ \sum_{k=1}^{2n-1} Z_{n,k}^{6V}=Z_n^{6V} $$
we have in the presence of a non-trivial $w$:
$$ \sum_{k=1}^{n} Z_{n,2k-1}^{6V} ({\bar b}_o{\bar b}_e)^{k-1} {\bar c}_o ({\bar a}_e{\bar a}_o)^{n-k} +\sum_{k=1}^{n-1} Z_{n,2k}^{6V} ({\bar b}_o{\bar b}_e)^{k-1} {\bar b}_o {\bar c}_e {\bar a}_o ({\bar a}_e{\bar a}_o)^{n-k-1}=Z_n^{6V}[w]$$
Noting that ${\bar c}_e={\bar c}_o$, ${\bar b}_o{\bar a}_o={\bar a}_e{\bar a}_o={\bar \omega}_0$ and ${\bar b}_o{\bar b}_e={\bar \omega}_1$, and using the parameter $\tau$ of \eqref{deftau}, we finally get an expression for the generating function of the sums of refined partition functions 
$Z_{n,2k-1}^{6V}+Z_{n,2k}^{6V}$:
\begin{equation}\label{sumpf}Z_n^{6V}(\tau):= \sum_{k=1}^{n} (Z_{n,2k-1}^{6V}+Z_{n,2k}^{6V})\, \tau^{k-1}=\frac{1}{{\bar \omega_0}^{n-1} \sqrt{-w}}\, Z_n^{6V}[w]
\end{equation}
or equivalently the generating function $h_n^{6V}(\tau)$ for sums of one-point functions $ \frac{Z_{n,2k-1}^{6V}+Z_{n,2k}^{6V}}{Z_n^{6V}}$:
\begin{equation}\label{sumrule3}h_n^{6V}(\tau):=\sum_{k=1}^{n} \frac{Z_{n,2k-1}^{6V}+Z_{n,2k}^{6V}}{Z_n^{6V}} \,\tau^{k-1}
=\frac{Z_n^{6V}(\tau)}{Z_n^{6V}}=\frac{1}{{\bar \omega_0}^{n-1} \sqrt{-w}}\, \frac{Z_n^{6V}[w]}{Z_n^{6V}}
\end{equation}
where we set $Z_{n,2n}^{6V}=0$ by convention.

Combining \eqref{sumrule2}, \eqref{sumrule3} and the relation \eqref{rela206}, we finally get the following relation between one-point generating functions of the $20V$ and $6V$ models:
\begin{equation}\label{relasixt}
h_{n}^{20V}(\tau)={\bar \omega}_0^{n-1}\,\frac{q^2-q^{-2}w}{\sqrt{2}}\left(\frac{1-w}{2\,{\bar \omega}_0}\right)^{2n-1} \, h_{n}^{6V}(\tau)
=\left(\frac{1+\tau}{2}\right)^{n-1}  \,h_{n}^{6V}(\tau)
\end{equation}
where we have identified the quantity $\frac{1-w}{\sqrt{2}(q^2-q^{-2}w)}=\frac{1+\tau}{2}$.

\subsection{Relation to the inhomogeneous U-turn boundary partition function}

Let us finally relate the above quantities to the U-turn boundary partition function $Z_n^{6V-U}[w]$.
First we note that \eqref{rela6VU} implies:
$$\frac{Z_{n}^{6V}[w]}{Z_{n}^{6V}}=\frac{\sqrt{2}}{q^{-2}-q^2w}\,  \frac{Z_{n}^{6V-U}[w]}{Z_{n}^{6V-U}} $$
and moreover that \eqref{Udelta} implies the relation:
$$\frac{Z_{n}^{6V-U}[w]}{Z_{n}^{6V-U}}=\sqrt{-w}\left(\frac{1-w}{2}\right)^{2n}\left( \frac{q^{-2}-q^2 w}{\sqrt{2}}\right)^{n+1} \left( \frac{q^{2}-q^{-2} w}{\sqrt{2}}\right)^{n} \frac{\Delta_n[w]}{\Delta_n}$$
where we denote for short $\Delta_n[w]$ the quantity $\Delta_n[z_1,...,z_n;w_1,...,w_n]$ \eqref{defdelt} with all $z_i=q^2$ and all $w_i=-1$ except $w_n=w$, and $\Delta_n:=\Delta_n[-1]$.
Consequently, we have:
\begin{eqnarray}\label{sixvone}
h_n^{6V}(\tau)&=& \frac{1}{{\bar \omega_0}^{n-1} \sqrt{-w}}
\frac{\sqrt{2}}{q^{-2}-q^2w}\,  \frac{Z_{n}^{6V-U}[w]}{Z_{2n\times n}^{6V-U}} \nonumber \\
&=&\left(\frac{1-w}{2}\right)^{n+1}\left(\frac{q^2-q^{-2}w}{\sqrt{2}}\right)\left(\frac{q^{-2}-q^{2}w}{\sqrt{2}}\right)^n 
\frac{\Delta_n[w]}{\Delta_n}
\end{eqnarray}
Similarly:
\begin{equation}\label{20vonept}
h_n^{20V}(\tau)=\left(\frac{1-w}{2}\right)^{2n}\left(\frac{q^2-q^{-2}w}{\sqrt{2}}\right)^{2}\left(\frac{q^{-2}-q^{2}w}{q^2-q^{-2}w}\right)^n \frac{\Delta_n[w]}{\Delta_n}
\end{equation}

\subsection{Refined U-turn boundary}
We now turn to the calculation of the refined determinant $\Delta_n[w]$. We start with a refinement of Theorem \ref{thmdelta}.
We consider the semi-homogeneous limit
$$\Delta_n'[z,w]:=\lim_{z_1,...,z_n\to z\atop w_1,...,w_{n-1}\to -1, w_n\to w}\Delta_n[z_1,...,z_n;w_1,...,w_n]$$
As in the proof of Theorem \ref{thmdelta}, we proceed in two steps. First, we take the successive limits $w_1,w_2,...,w_{n-1}\to -1$.
This step is rigorously identical, as only the limit $w_n\to w$ is different. Let us consider:
\begin{eqnarray*}
\Delta_n'[\bz,w]&:=&\lim_{w_1,...,w_{n-1}\to -1\atop  w_n\to w}\Delta_n[z_1,...,z_n;w_1,...,w_n]\\
&=&\lim_{v_1,v_2,...,v_n\to 0}\frac{\det\left( \bF[v_1]\bF[v_2] \cdots \bF[v_{n-1}] {\tilde \bF}[v_n]\right) }{2^{n-1}(w+1)^{2n-2}(1-w^2)\prod_{1\leq i<j\leq n-1} (v_j^2-v_i^2) \prod_{i=1}^{n-1} v_i }
\end{eqnarray*}
where $\bF[v]$ is as in Sect.~\ref{homsec}, and where ${\tilde \bF}[v]:=\bF[v+1+w]$, which accounts for taking the limit $w_n\to w$ instead of $-1$ in the last column. Writing $\bF[v]=v\bG[v]$ as before, and performing successively the limits $v_1\to 0, ...,v_{n-1}\to 0$,
we end up with
$$\Delta_n'[\bz,w]=\frac{\det\left( \bF_1\bF_2 \cdots \bF_{n-1} {\bF}[w+1]\right) }{2^{n-1}(w+1)^{2n-1}(1-w)}$$

The second step consists in computing:
\begin{eqnarray*}
\Delta_n'[z,w]&=&\lim_{z_1,...,z_n\to z} \Delta_n'[\bz,w]\\
&=&\lim_{u_1,...,u_n\to 0}
 \frac{\det\left( \bF_1\bF_2 \cdots \bF_{n-1} {\bF}[w+1]\right) }
 {2^{n-1}(w+1)^{2n-1}(1-w)(1-z^2)^{n(n+1)/2} \prod_{1\leq i<j \leq n} (u_i-u_j)} 
\end{eqnarray*}
Like in Sect. \ref{homsec}, let us denote by ${\tilde \bR}[u_1], ...,{\tilde \bR}[u_n]$ the transpose of the rows of the matrix 
$\left( \bF_1\bF_2 \cdots \bF_{n-1} {\bF}[w+1]\right)$, where ${\tilde \bR}[u]$ has entries 
$(f_1[z+u], f_3[z+u],...,f_{2n-3}[z+u],f[z+u;w+1])$,
with $f[z+u;w+1]=f_U(u,w+1)$ \eqref{fu} and $f_i$ as in \eqref{columns}. Expanding again ${\tilde \bR}[u]=\sum_{m\geq 0}u^m\,{\tilde \bR}_m$, and performing the successive limits $u_1\to 0, ...,u_n\to 0$, we arrive at:
$$\Delta_n'[z,w]= (-1)^{n(n-1)/2}\, \frac{\det\left( ({\tilde \bR}_1{\tilde \bR}_2 \cdots {\tilde \bR}_{n})^t\right) }
 {2^{n-1}(w+1)^{2n-1}(1-w)(1-z^2)^{n(n+1)/2}} 
$$
The entries $( {\tilde \bR}_m)_i$ of the vector $ {\tilde \bR}_m$ are easily identified as 
$f_U(u,v)\vert_{u^iv^{2m-1}}$ for $m=1,2,...,n-1$ and
$f_U(u,w+1)\vert_{u^i}$ for $m=n$.
This result may be recast into the following.

\begin{lemma}
We have
$$ \Delta_n'[z,w]= (-1)^{n(n-1)/2}\, \frac{\det_{1\leq i,j\leq n-1}\left( 
f_U(z,w;u,v)\vert_{u^iv^{2j+1}}\right) }
 {2^{n-1}(w+1)^{2n-1}(1-w)(1-z^2)^{n(n+1)/2}} $$
in terms of the function:
\begin{equation}\label{genoU}
f_{U}(z,w;u,v)
:= f_U(z;u,v)+v^{2n-1} \left(f_U(z;u,w+1) -f_U(z;u,v)\vert_{v^{2n-1}}\right)   
\end{equation}
\end{lemma}
\noindent (Note the subtraction of the contribution of order $v^{2n-1}$ of $f_U(z;u,v)$, to avoid over-counting.).
We are now ready to compute the quantity $\Delta_n[w]=\Delta_n'[q^2,w]$ corresponding to $z=q^2$. 
\begin{thm}\label{finthmdelta}
We have the following determinant formula:
\begin{eqnarray}
\quad\qquad \Delta_n[w] &=& \frac{q^{6n(n+1)}}{2^{\frac{n(7n+5)}{4}}} \left(\frac{2}{1-w}\right)^{2n} \, \det_{0\leq i,j\leq n-1}\left( 
g_U'(u,v)\vert_{u^iv^j} \right) \label{deltafin}\\
\quad\qquad g_U'(u,v)&=& g_U(u,v)+v^{n-1}\left(\frac{1+u}{1-u}\right)^{2n} \left\{ \frac{4(1-u)^2}{(1-u)^2(1-w)^2+(1+u)^2(1+w)^2}-1\right\}
\label{guprime}
\end{eqnarray}
with $g_U(u,v)$ as in \eqref{gu}.
\end{thm}
\begin{proof}
Like in Sect.~\ref{homsec}, we apply Lemma \ref{lemdeltatwo} to the function $f_U(q^2,w;u,v)$ \eqref{genoU}, with the same choices of $a,b,\al,\beta,\gamma,\delta$. More precisely we apply the substitutions:
\begin{equation}\label{substi}
u\to \frac{2q^6  u}{1-q^4 u},\qquad v\to \frac{2 q^{-4} v}{1+q^{-4}v}, \qquad w= \frac{2 q^{-4} x}{1+q^{-4}x}-1
\end{equation}
and we multiply it by an overall factor of $\frac{1+q^4 u}{(1-q^4 u)(1+q^{-4}v)^2} $. The last substitution
in \eqref{substi} is simply a change of variable $w\to x$, which allows to treat $v$ and $w+1$ similarly. 
We first compute the effect of this transformation on the term $v^{2n-1} f_U(z;u,w+1)$, up to higher orders in $v$:
\begin{eqnarray*}
&&\frac{1+q^4 u}{(1-q^4 u)(1+q^{-4}v)^2} (2 q^{-4}v)^{2n-1} f_U\left(q^2;\frac{2q^6  u}{1-q^4 u},\frac{2 q^{-4} x}{1+q^{-4}x}\right)\\
&&\quad = (2v)^{2n-1}q^4(-1)^{n}\frac{q^4}{1-q^4}\, (1+q^{-4}x)^2\, x\, g_U(u,x^2)+O(v^{2n})
\end{eqnarray*}
with $g_U(u,v)$ as in \eqref{gu}.
The complete function $f_{U}(q^2,w;u,v)$ is therefore transformed into:
\begin{eqnarray*}
{\tilde f}_U(q^2,x;u,v)
&=& v\frac{q^4}{1-q^4} \left\{g_U(u,v^2) \right.\\
&&\qquad \left.+ v^{2n-2} \left(2^{2n-1} q^4(-1)^{n}(1+q^{-4}x)^2 x\, g_U(u,x^2) -g_U(u,v^2)\vert_{v^{2n-2}}\right)\right\}
\end{eqnarray*}
In this last expression we may replace $x g_U(u,x^2)=\sum_{m\geq 0} x^{2m+1} 
\left\{\left(\frac{1+u}{1-u}\right)^{2m+2}-u^{2m+2}\right\}$ with the truncated expansion
$$x {\bar g}(u,x^2)=\sum_{m\geq n-1} x^{2m+1} \left\{\left(\frac{1+u}{1-u}\right)^{2m+2}-u^{2m+2}\right\} 
=x^{2n-1}\left(\frac{1+u}{1-u}\right)^{2n}\, \frac{1}{1-x^2 \left(\frac{1+u}{1-u}\right)^{2}}  +O(u^{2n})$$
Indeed, the first $n-1$ terms in the series produce columns that are linear combinations of the first $n-1$ columns in the determinant,
and can therefore be omitted without affecting the result. Finally we restore the variable $w$ by inverting the change of variables
in \eqref{substi} into $x=q^4\frac{1+w}{1-w}$. The corresponding function finally reads:
\begin{eqnarray*}
&&{\tilde f}_U'(q^2,w;u,v)=-\frac{v q^{-2}}{\sqrt{2}} \left\{g_U(u,v^2)\right. \\
&&\quad \left. +v^{2n-2} \left(\frac{1+u}{1-u}\right)^{2n}\left(\left(\frac{2(1+w)}{1-w}\right)^{2n-1}  \frac{4(1-u)^2}{(1-u)^2(1-w)^2+(1+u)^2(1+w)^2}-1\right)\right\} \\
&&\quad =:  -\frac{v q^{-2}}{\sqrt{2}}  {\tilde g}_U'(u,v^2) 
\end{eqnarray*}
By Lemma \ref{lemdeltatwo} we deduce:
\begin{eqnarray*}
\Delta_n[w] &=&\frac{q^{6n(n+1)}}{2^{\frac{n(7n+5)}{4}} \frac{1-w}{2} (1+w)^{2n-1}} \, \det_{0\leq i,j\leq n-1} \left(
{\tilde g}_U'(u,v)\vert_{u^iv^j}\right)\\
&=& \frac{q^{6n(n+1)}}{2^{\frac{n(7n+5)}{4}}} \left(\frac{2}{1-w}\right)^{2n} \, \det_{0\leq i,j\leq n-1} \left(
{ g}_U'(u,v)\vert_{u^iv^j}\right)
\end{eqnarray*}
in which we have divided explicitly the last column by $\left(\frac{2(1+w)}{1-w}\right)^{2n-1}$, thus allowing to replace ${\tilde g}_U'(u,v)$ with ${ g}_U'(u,v)$.
The theorem follows. 
\end{proof}

\subsection{Refined partition functions for 20V and 6V models}

Theorem \ref{finthmdelta} provides the following formula for the ratio $\Delta_n[w]/\Delta_n$:
$$\frac{\Delta_n[w]}{\Delta_n}=\Big(\frac{2}{1-w}\Big)^{2n}\, \, \frac{\det_{0\leq i,j\leq n-1} \left(
g_U'(u,v)\vert_{u^iv^j}\right)}{\det_{0\leq i,j\leq n-1} \left(
g_U(u,v)\vert_{u^iv^j}\right)}$$
in terms of the functions $g_U(u,v)$ \eqref{gu} and $g_U'(u,v)$ \eqref{guprime}.

Using the relations \eqref{sixvone} and \eqref{20vonept} and the result of Theorem \ref{finthmdelta}, we immediately get the following expressions for
the one-point function generating functions of the 6V and 20V models, in terms of the variable $\tau=\frac{1-q^4 w}{q^4-w}$ 
(i.e. $w=\frac{1-q^4 \tau}{q^4-\tau}$), and of the function $g_{20V}=g_U$:

\begin{thm}\label{oneptthm}
The one-point generating functions for the $6V$ and $20V$ models read  for all $n\geq 1$:
\begin{eqnarray}
h_n^{6V}(\tau)&=& \frac{\det_{0\leq i,j\leq n-1} \left(
g_{6V}^{\rm ref}(u,v)\vert_{u^iv^j}\right)}{\det_{0\leq i,j\leq n-1} \left(
g_{20V}(u,v)\vert_{u^iv^j}\right)} \label{sixVonept}\\
h_n^{20V}(\tau)&=&\frac{\det_{0\leq i,j\leq n-1} \left(
g_{20V}^{\rm ref}(u,v)\vert_{u^iv^j}\right)}{\det_{0\leq i,j\leq n-1} \left(
g_{20V}(u,v)\vert_{u^iv^j}\right)}\label{twentyVonept}
\end{eqnarray}
in terms of the functions:
\begin{eqnarray}
\qquad g_{6V}^{\rm ref}(u,v)&=& g_{20V}(u,v)+v^{n-1}  \left(\frac{1+u}{1-u}\right)^{2n} \left\{ \left(\frac{2}{1+\tau}\right)^{n-1} \frac{\tau^n\, (1-u)^2}{(\tau-u)(1-\tau u)} -1\right\} \label{6vgfone}\\
\qquad g_{20V}^{\rm ref}(u,v)&=& g_{20V}(u,v)+v^{n-1}  \left(\frac{1+u}{1-u}\right)^{2n} \left\{\frac{\tau^n\, (1-u)^2}{(\tau-u)(1-\tau u)} -1\right\} 
\label{twentyvgfone}
\end{eqnarray}
with $g_{20V}(u,v)$ as in \eqref{f20v}.
\end{thm}

In particular, using Theorem \ref{thm20V} to identify the denominator of \eqref{twentyVonept} as $Z_n^{20V}$,  this implies immediately the following:
\begin{thm}\label{refcor}
The refined partition function generating function for the $20V$ model reads:
$$Z_n^{20V}(\tau)=\det_{0\leq i,j\leq n-1} \left(
g_{20V}^{\rm ref}(u,v)\vert_{u^iv^j}\right)$$
with $g_{20V}^{\rm ref}$ as in \eqref{twentyvgfone}.
\end{thm}
\begin{remark}\label{ref20vrem}
We may interpret again the result as:
\begin{equation}\label{inter20v}
Z_n^{20V}(\tau)=\det\left( (P^{(n)}(\tau))_n\right)  
\end{equation} 
the determinant of a finite truncation of an infinite matrix $P^{(n)}(\tau)$, which differs from the matrix $P$ of Remark \ref{rem20V} only in its $n$-th column, and is generated by $f_{P^{(n)}(\tau)}(u,v)=g_{20V}^{\rm ref}(u,v)$ \eqref{twentyvgfone}.
\end{remark}

Direct calculation using Theorem \ref{refcor} gives the following first few values for $n=1,2...,5$
of the refined partition function generating polynomials $Z_{n}^{20V}(\tau)$
for the 20V model on the quadrangle $\cQ_n$:
\begin{eqnarray*}
Z_1^{20V}(\tau)&=&1\\
Z_2^{20V}(\tau)&=&1 + 2\tau +\tau^2\\
Z_3^{20V}(\tau)&=&4 + 15 \tau + 22 \tau^2 +15 \tau^3 + 4 \tau^4\\
Z_4^{20V}(\tau)&=&60 + 328\tau+ 772\tau^2 + 1008 \tau^3 + 772 \tau^4 + 328 \tau^5 + 60 \tau^6\\
Z_5^{20V}(\tau)&=&3328 + 23868 \tau+ 76856\tau^2 + 145860\tau^3 + 179088 \tau^4 + 
145860\tau^5 \\
&&+ 76856 \tau^6 + 23868 \tau^7 + 3328\tau^8
\end{eqnarray*}
\begin{remark}\label{palinrem}
Note that $Z_n^{20V}(\tau)$ is a palindromic polynomial of degree $2n-2$, namely that $Z_{n,2n-k}^{20V}=Z_{n,k}^{20V}$ for all $k=1,2,...,2n-1$. 
Indeed, from \eqref{twentyvgfone} we see that the only dependence on 
$\tau$ comes from the $n$-th column of $P^{(n)}(\tau)$, generated by $\varphi_n(\tau):=\left(\frac{1+u}{1-u}\right)^{2n} \frac{\tau^n\, (1-u)^2}{(\tau-u)(1-\tau u)}$. 
The property follows from  the fact that $\varphi_n(\tau^{-1})=\tau^{-2(n-1)}\varphi_n(\tau)$.

Observe also that
\begin{equation}\label{obs}
Z_n^{20V}(0)=Z_{n,1}^{20V}=Z_{n,2n-1}^{20V}=Z_{n-1}^{20V} 
\end{equation}
Indeed,  the contributions to $Z_{n,2n-1}^{20V}$ all have a trivial topmost path made of $n$ horizontal steps in the top row, followed by $2n-1$ down steps along the rightmost vertical until its endpoint. This effectively cuts out the top row and right column from the domain $\cQ_n$, and the remaining paths form a configuration on $\cQ_{n-1}$. Similarly, the contributions to $Z_{n,1}^{20V}=Z_n^{20V}(0)$ all have the bottom-most diagonal entirely occupied by paths, while the rightmost vertical is empty except for the bottom-most edge. This effectively cuts out the bottom diagonal and right vertical lines of vertices, leaving us with paths on $\cQ_{n-1}$ as well.
\end{remark}

Finally, Theorem \ref{oneptthm} gives the corresponding $6V$ one-point functions:
\begin{eqnarray*}
h_1^{6V}(\tau)&=&1\\
h_2^{6V}(\tau)&=&\frac{1 + \tau }{2}\\
h_3^{6V}(\tau)&=&\frac{4 + 7 \tau + 4 \tau^2}{15}\\
h_4^{6V}(\tau)&=&\frac{15 + 37\tau+ 37\tau^2 + 15 \tau^3}{104}\\
h_5^{6V}(\tau)&=&\frac{64 + 
 203 \tau+ 282\tau^2 + 203\tau^3 + 64 \tau^4}{816}
\end{eqnarray*}
We note again that these polynomials are palindromic, a property inherited from the 20V one-point functions via the relation \eqref{relasixt}. In particular, comparing the
top and constant coefficients we get a non-trivial relation: $Z_{n,1}^{6V}+Z_{n,2}^{6V}=Z_{n,2n-1}^{6V}$.

\section{Refined domino tilings of the Aztec triangle and connection to 20V refined partition functions}
\label{refdtsec}

\subsection{Refined domino tiling partition functions}

\begin{figure}
\begin{center}
\includegraphics[width=9cm]{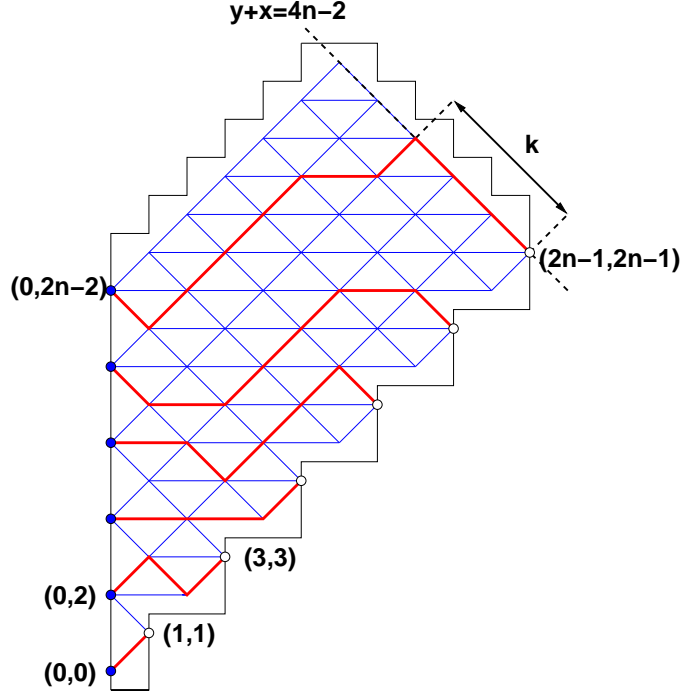}
\end{center}
\caption{\small A typical configuration of non-intersecting Schr\"oder paths for the domino tiling model of the Aztec triangle $\cT_n$, contributing to the refined partition function $Z_{n,k}^{DT}$, i.e. such that the first entry point of any path into the line $x+y=4n-2$ lies at some distance $k$ of the endpoint $(2n-1,2n-1)$.}
\label{fig:pachteref}
\end{figure}

Consider again the path model for the domino tilings of the Aztec triangle ${\mathcal T}_n$. In any non-intersecting path configuration, the path ending at the topmost point $(2n-1,2n-1)$ must first hit the diagonal line $x+y=4n-2$ at some point $(2n-1-k,2n-1+k)$ for some $k\in [0,n-1]$. Let us denote by $Z^{DT}_{n,k}$ the corresponding contribution to the total partition function (see Fig.~\ref{fig:pachteref} for an illustration).

Alternatively, let us denote by $W^{DT}_{n,k}$ the partition function for
non-intersecting Schr\"oder paths with starting points $(0,2i)$, $i=0,1,...,n-1$ and endpoints $(2j+1,2j+1)$, $j=0,1,...,n-2$ together with
the point $(2n-1-k,2n-1+k)$. 
The condition of first hitting the line $y+x=4n-2$ at $(2n-1-k,2n-1+k)$ is equivalent to imposing that the path visits the point $(2n-1-k,2n-1+k)$ but does not visit the point $(2n-1-k-1,2n-1+k+1)$. Note that once a point is visited on the line $x+y=4n-2$, the rest of the path descends along this line until it reaches the endpoint $(2n-1,2n-1)$. This implies immediately for all $k\in [0,n-1]$:
$$Z^{DT}_{n,k} =W^{DT}_{n,k}-W^{DT}_{n,k+1}$$
For convenience, we define the generating function:
$$Z^{DT}(t):=\sum_{k=0}^{n-1} Z^{DT}_{n,k} \, t^k$$
for refined domino tiling partition functions of the Aztec triangle $\cT_n$.
We have the following:

\begin{thm}\label{refdominothm}
The generating function $Z^{DT}(t)$ reads:
\begin{eqnarray}Z^{DT}(t)&=&\det_{0\leq i,j\leq n-1} \left(
f_{DT}(t;u,v)\Big\vert_{u^iv^j}\right) \label{DTref} \\
f_{DT}^{\rm ref}(u,v)&=&f_{DT}(u,v)+v^{n-1} \frac{\al_+(u)^{2n}}{\sqrt{1+6u+u^2}}\, \frac{(t-1)u}{\al_+(u)-t u}\label{refat}
\end{eqnarray}
where $\al_{+}(u)= \frac{1}{2}\big(1 + u+\sqrt{1 + 6 u + u^2}\big)$ and $f_{DT}(u,v)$ as in \eqref{fat}.
\end{thm}
\begin{proof}
Let us first compute $W^{DT}_{n,k}$.
By direct application of the Lindstr\"om-Gessel-Viennot theorem \cite{LGV1,GV}, we have $W^{DT}_{n,k}=\det(M_n^{(k)})$, where the $n\times n$ matrix 
$M_n^{(k)}$ differs from the matrix $M_n$ (of Lemma \ref{thmzat}, Remark \ref{remAT} and  Theorem \ref{genDT}) only in its
last column, where the entry $M_{i,n-1}=\Sigma_{2n-1,2n-1-2i}$ \eqref{sigentry} must be replaced by $(M_n^{(k)})_{i,n-1}=\Sigma_{2n-1-k,2n-1-2i+k}$ as the corresponding endpoint has been shifted by $(-k,k)$.
Recall that $\Sigma_{i,j}=\sigma(u,v)\vert_{u^{\frac{i+j}{2}}v^{\frac{i-j}{2}}}$, with $\sigma(u,v)$ as in \eqref{defsig}.
We deduce that
\begin{eqnarray*}\Sigma_{2n-1-k,2n-1-2i+k}&=&\oint \frac{dx}{2i\pi x} \frac{dy}{2i\pi y} \frac{1}{x^{2n-1-i}y^{i-k}}\frac{1}{1-x-y-x y}\\
&=&
\oint \frac{du}{2i\pi u} \frac{dv}{2i\pi v} \frac{1}{v^{2n-1-k}u^{i-k}}\frac{1}{1-v(1+u)-u v^2}
\end{eqnarray*}
by changing again to variables $u=y/x$, $v=x$. Noting that:
$$\frac{1}{1-v(1+u)-u v^2} =\frac{1}{v\sqrt{1+6u+u^2}}\left(\frac{1}{1-v \al_+(u)}-\frac{1}{1-v \al_-(u)} \right)$$
where
$$\al_{\pm}(u)= \frac{1}{2}\big(1 + u\pm \sqrt{1 + 6 u + u^2}\big)$$
we get
$$\frac{1}{1-v(1+u)-u v^2}\Bigg\vert_{v^{2n-1-k}}=\frac{\al_+(u)^{2n-k}-\al_-(u)^{2n-k} }{\sqrt{1+6u+u^2}} $$
Note that $\al_-(u)=-u+O(u^2)$ hence the term $\al_-(u)^{2n-k}=O(u^{2n-k})$ does not contribute to the coefficient of $u^{i-k}$.
The generating function for the $n$-th column of $M_n^{(k)}$ reads:
$$v^{n-1}\sum_{i=0}^{n-1} \Sigma_{2n-1-k,2n-1-2i+k} \, u^i =v^{n-1} u^k \frac{\al_+(u)^{2n-k}}{\sqrt{1+6u+u^2}} +O(u^n)$$
The matrix $M_n^{(k)}$ is the finite truncation to the first $n$ rows and columns of the infinite matrix $M^{(n,k)}$,
with generating function
$$f_{M^{(n,k)}}(u,v)= f_M(u,v)+v^{n-1} \frac{u^k \,\al_+(u)^{2n-k}-\al_+(u)^{2n}}{\sqrt{1+6u+u^2}} $$

Using the linearity of the determinant w.r.t. its last column, we immediately deduce that
$Z_{n,k}^{DT}=W^{DT}_{n,k}-W^{DT}_{n,k+1}=\det({\tilde M}_n^{(k)})$ where the $n\times n$ matrix ${\tilde M}_n^{(k)}:=({\tilde M}^{(n,k)})_n$ is the finite truncation of the infinite matrix
${\tilde M}^{(n,k)}$ with generating function:
$$f_{{\tilde M}^{(n,k)}}(u,v)= f_M(u,v)+v^{n-1} \frac{u^k \al_+(u)^{2n-k}-u^{k+1} \al_+(u)^{2n-k-1}-\al_+(u)^{2n}}{\sqrt{1+6u+u^2}} $$
Similarly, the generating function $Z^{DT}(t)=\det\Big(({\tilde M}^{(n)}(t))_n\Big)$ where the $n\times n$ matrix $({\tilde M}^{(n)}(t))_n$
is the finite truncation of the infinite matrix ${\tilde M}^{(n)}(t)$ generated by:
\begin{eqnarray*}f_{{\tilde M}^{(n)}(t)}(u,v)&=&f_M(u,v)+v^{n-1}  \frac{\sum_{k\geq 0}t^k(u^k \al_+(u)^{2n-k}-u^{k+1} \al_+(u)^{2n-k-1})-\al_+(u)^{2n}}{\sqrt{1+6u+u^2}}\\
&=&f_M(u,v)+v^{n-1} \frac{\al_+(u)^{2n}}{\sqrt{1+6u+u^2}}\, \frac{(t-1)u}{\al_+(u)-t u}=f_{DT}^{\rm ref}(u,v)
\end{eqnarray*}
(Here we have considered the sum over all $k\geq 0$, as the terms of order $u^n$ and above do not affect the truncation.)
This completes the proof of Theorem \ref{refdominothm}.
\end{proof}

Direct calculation using Theorem \ref{refdominothm} leads to the following first few values of the refined partition function generating functions $Z_n^{DT}(t)$, up to $n=5$:
\begin{eqnarray*}
Z_1^{DT}(t)&=&1\\
Z_2^{DT}(t)&=&3+t\\
Z_3^{DT}(t)&=&37 + 19 t + 4 t^2\\
Z_4^{DT}(t)&=&1780 + 1100 t + 388 t^2 + 60 t^3\\
Z_5^{DT}(t)&=&324948 + 222716 t + 100724 t^2 + 27196 t^3 + 3328 t^4
\end{eqnarray*}
We note the obvious sum rule $Z_n^{DT}(1)=Z_n^{DT}$.
Observe also that the leading coefficient of $t^{n-1}$ in $Z_n^{DT}(t)$ is nothing but 
\begin{equation}\label{minusone}Z_{n-1,1}^{DT}=Z_{n-1}^{DT}
\end{equation}
Indeed, imposing that the topmost path 
(from $(0,2n-2)$ to $(2n-1,2n-1)$) visits the topmost vertex on the line $x+y=4n-2$ forces this path to be a trivial succession of $n$ up steps, followed by $n-1$ down steps, effectively reducing the domain $\cT_n$ to $\cT_{n-1}$ for the $n-1$ remaining paths in each such configuration.

\subsection{Relation between the refined tiling and refined 20V partition functions}
\label{seclutfin}

\begin{thm}\label{fullequivthm}
We have the following relation between refined partition functions of the $20V$ model on the quadrangle $\cQ_n$ and of the domino tilings of the Aztec triangle $\cT_n$:
\begin{equation}\label{relaraf} Z_n^{20V}(\tau)=\tau^n\,  \frac{ Z^{DT}(\tau)+\tau^{-1}Z^{DT}(\tau^{-1})}{1+\tau} \end{equation}
\end{thm}
\begin{proof}
We apply the same technique as in the case of ordinary partition functions. Let $L$ be the infinite matrix \eqref{defL} of Sect. \ref{equivsec} and ${\tilde M}^{(n)}(t)$ generated by $f_{{\tilde M}^{(n)}(t)}(u,v)=f_{DT}^{\rm ref}(u,v)$ of Theorem \ref{refdominothm}. We compute:
\begin{eqnarray}
\label{wecompute}
f_{L{\tilde M}^{(n)}(t)}(u,v)&=& \frac{1+2u-u^2}{1-u}f_{{\tilde M}^{(n)}(t)}\left(u\frac{1+u}{1-u},v\right)=f_P(u,v)+v^{n-1} \left(\frac{1+u}{1-u}\right)^{2n}\frac{(t-1)u}{1-t u}\nonumber \\
&=&g_{20V}(u,v)+v^{n-1} \left(\frac{1+u}{1-u}\right)^{2n}\left\{ \frac{1-u}{1-t u} -1\right\}
\end{eqnarray}
by use of the result $LM=P$ of Section \ref{equivsec}.
Let us now denote by $\bC(t)$ the $n$-th column of the matrix $L{\tilde M}^{(n)}(t)$. We have $Z^{DT}(t)=\det\left((L{\tilde M}^{(n)}(t))_n\right)$.
By linearity of the determinant w.r.t. the $n$-th column, the quantity $\frac{t^n}{1+t}\big( Z^{DT}(t)+t^{-1}Z^{DT}(t^{-1})\big)$ 
is the determinant of the finite truncation of
an infinite matrix ${\bar  M}^{(n)}(t)$ defined as follows: it is identical to $L{\tilde M}^{(n)}(t)$ except for the $n$-th column, where $\bC(t)$ is replaced with the combination $\frac{t^n}{1+t} \{ \bC(t)+t^{-1}\bC(t^{-1})\}$.
Consequently:
\begin{eqnarray*}t^n\, \frac{ Z^{DT}(t)+t^{-1}Z^{DT}(t^{-1})}{t+1}&=&\det\left( ({\bar M}^{(n)}(t))_n \right)\\
f_{{\bar M}^{(n)}(t)}(u,v)&=&g_{20V}(u,v)+v^{n-1} \left(\frac{1+u}{1-u}\right)^{2n} \left\{ t^n\,\frac{ \frac{1-u}{1-t u}+ \frac{1-u}{t-u}  }{1+t}-1\right\}\\
&=&g_{20V}^{\rm ref}(u,v)
\end{eqnarray*}
where we have recognized the generating function $g^{\rm ref}_{20V}(u,v)$ (at $\tau=t$) of Theorem \ref{oneptthm}.
The theorem follows by setting $t=\tau$, and using the identification ${\bar M}^{(n)}(\tau)=P^{(n)}(\tau)$ given by \eqref{inter20v}.
\end{proof}

\begin{remark}\label{reciprem}
Note that the result above gives an independent confirmation that the polynomial $Z_n^{20V}(\tau)$ is palindromic, 
namely that $\tau^{2n-2}Z_n^{20V}(\tau^{-1})=Z_n^{20V}(\tau)$ as noted in Remark \ref{palinrem}.
\end{remark}

Finally Eqn.\eqref{relaraf} implies the following relations between the refined partition functions of the $20V$ model and the domino tiling problem.
\begin{cor}\label{ref20dt}
The refined partition functions for the $20V$ model and for the domino tilings of the Aztec triangle are related as follows:
\begin{equation}\label{usefuldt} Z_{n,k}^{DT}= Z_{n,n+k+1}^{20V}+Z_{n,n+k}^{20V}=Z_{n,n-k-1}^{20V}+Z_{n,n-k}^{20V}
\end{equation}
for $k=0,1,...,n-1$ with the convention that $Z_{n,0}^{20V}=Z_{n,2n}^{20V}=0$.
\end{cor}

\section{Pentagonal 20V and truncated Aztec triangle domino tilings}
\label{conjsec}

In this section, we prove the correspondence (20V DWBC3-DT conjecture of \cite{DG19}) between 20V configurations on the domains
$\cP_{n,n-k}$ (see Fig.~\ref{fig:pachtri}) and the domino tiling configurations of the domains $\cT_{n,n-k}$ for $k=2$ and $3$. The method of proof is purely combinatorial and does not rely on integrability, except for its use of the $k=1$ results.

\subsection{Proof of the 20V-DT conjecture for $k=2$}

\begin{figure}
\begin{center}
\includegraphics[width=15cm]{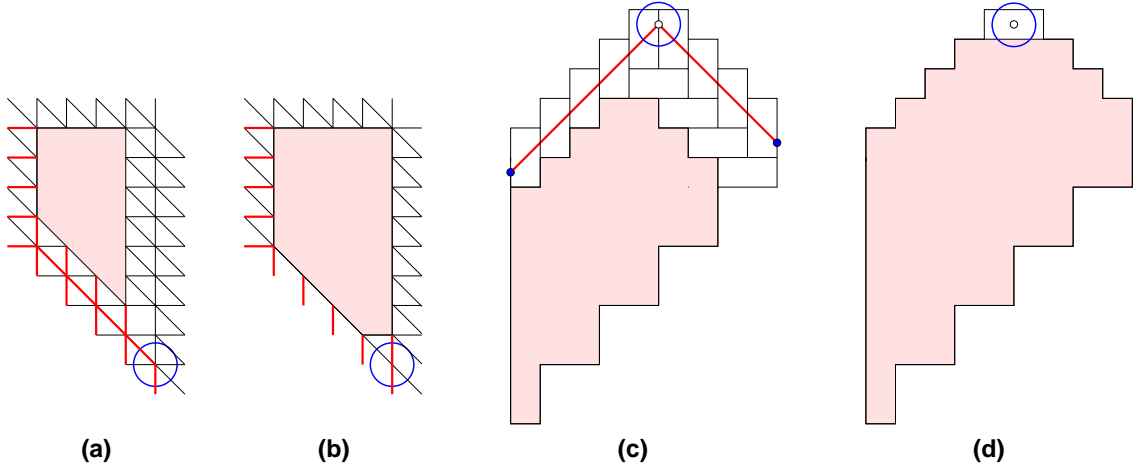}
\end{center}
\caption{\small The proof of the 20V-DT conjecture for $k=2$. The 20V partition function on $\cQ_n$ is decomposed according to the two possible configurations of the bottom vertex (a) and (b). The domino tiling partition function on $\cT_n$ is similarly decomposed according to whether the topmost vertex is visited (c) or not (d).}
\label{fig:pachter1}
\end{figure}

\begin{thm}
The number of configurations of the 20V DWBC3 model on the pentagon $\cP_{n,n-2}$ is identical to the number of domino tilings of the truncated Aztec triangle $\cT_{n,n-2}$, namely:
$$ Z^{20V}({\mathcal P}_{n,n-2})=Z^{DT}({\mathcal T}_{n,n-2}) $$
\end{thm}
\begin{proof}
The correspondence is proved by expressing  both relevant partition function in terms of known objects. 
For the 20V model on $\cP_{n,n-2}$, let us start from the partition function $Z_n^{20V}$ on the larger domain $\cQ_n=\cP_{n,n-1}$.
This quantity splits into two parts, according to the two possible local configuration of the bottom-most vertex displayed in Figs.\ref{fig:pachter1} (a) and (b) (circled vertex). In the case (a), the configuration of the bottom vertex induces a chain of forced edge occupations, due to the SW boundary condition. The contribution of case (a) is readily seen to reduce to 
$Z_{n,1}^{20V}=Z_{n-1}^{20V}$. The case (b) contributes exactly $Z^{20V}({\mathcal P}_{n,n-2})$. We deduce the relation:
$$ Z^{20V}({\mathcal P}_{n,n-2})= Z^{20V}_n-Z^{20V}_{n-1}$$ 
For the domino tiling of $\cT_{n,n-2}$ let us also start from the partition function $Z_n^{DT}$ on the larger domain $\cT_n=\cT_{n,n-1}$ in the lattice path formulation. Analogously, this quantity splits into two parts according to whether the topmost vertex is visited or not, as displayed in Figs.\ref{fig:pachter1} (c) and (d). The case (c) corresponds to a unique configuration of the topmost path, and the contribution reduces to $Z_{n-1}^{DT}$. The case (d) contributes precisely $Z^{DT}({\mathcal T}_{n,n-2})$. This gives the relation
$$ Z^{DT}({\mathcal T}_{n,n-2})= Z^{DT}_n-Z^{DT}_{n-1}$$
The theorem follows from the result of Theorem \ref{equivalencethm}, i.e. the identities $Z^{20V}_m=Z^{DT}_m$ for all $m$.
\end{proof}

\subsection{Proof of the 20V-DT conjecture for  $k=3$}

\begin{figure}
\begin{center}
\includegraphics[width=15cm]{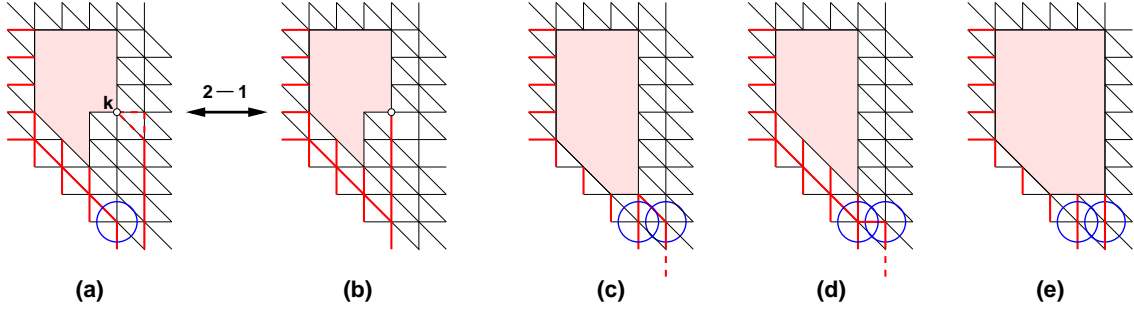}
\end{center}
\caption{\small Decomposition of $Z^{20V}({\mathcal P}_{n,n-2})$ into four parts (a), (c),(d), (e), according to the bottom vertex configurations (circled). Cases (a) are further decomposed according to the point of exit of the rightmost path at height $k$ (empty circle) from the second rightmost vertical line.
There are $2$ distinct configurations of the rightmost path after this exit point (paths using edges among the dashed ones), in $2$-1 correspondence with the partition function (b) on $\cQ_{n-1}$, equal to the sum of the refined partition function $\sum_{j\geq k} Z_{n,j}^{20V}$. Cases (c) and (d) sum up to 
$Z_{n,2}^{20V}-Z_{n,1}^{20V}$, while (e) contributes $Z^{20V}({\mathcal P}_{n,n-3})$.}
\label{fig:pachter2}
\end{figure}

\begin{figure}
\begin{center}
\includegraphics[width=14cm]{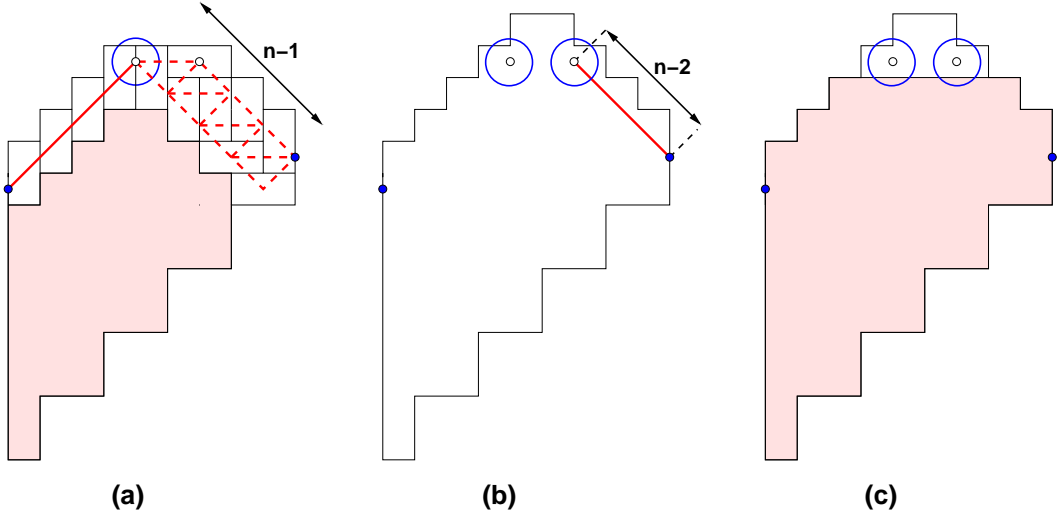}
\end{center}
\caption{\small Decomposition of $Z^{DT}({\mathcal P}_{n,n-2})$ into three pieces (a), (b), (c) according to the configurations of the two top vertices (empty circles). Case (a) gives rise to $2(n-1)$ times the domino problem on $\cT_{n-1}$, as there are $2(n-1)$ choices for how the top path ends, using only edges among the dashed ones. Case (b) receives contributions from the refined partition function where paths first enter the $x+y=4n-2$ line at distance $n-1$ from the top endpoint (i.e. at the right top vertex), but that do not visit the left vertex. Case (c) contributes $Z^{DT}({\mathcal T}_{n,n-3})$.}
\label{fig:pachter3}
\end{figure}

\begin{thm}
The number of configurations of the 20V DWBC3 model on the pentagon $\cP_{n,n-3}$ is identical to the number of domino tilings of the truncated Aztec triangle $\cT_{n,n-3}$, namely:
$$ Z^{20V}({\mathcal P}_{n,n-3})=Z^{DT}({\mathcal T}_{n,n-3})$$
\end{thm}
\begin{proof}
We proceed like in the case $k=2$. For the 20V model on $\cP_{n,n-3}$, let us start from the partition function 
$Z^{20V}({\mathcal P}_{n,n-2})$. It splits into four parts according to the local configuration of the two bottom vertices, as shown in Figs.\ref{fig:pachter2} (a), (c), (d) and  (e).

The first contribution (a) can be further split according to the configuration of the rightmost column: there are two ways in which the path can leave the second rightmost vertical line at an exit point say at height $k$ (by a horizontal or diagonal step, followed by vertical steps until the endpoint). Removing the last column, and inserting vertical steps up to the point of exit at height $k$ into the second vertical, we get a 2 to 1 mapping to refined configurations of the 20V on $\cQ_{n-1}$ that reach the last column at or above height $k$, resulting in a contribution $\sum_{j\geq k}Z_{n-1,j}^{20V}$ for each $k$, hence a total of $2\sum_{k=1}^{2n-3}\sum_{j\geq k}Z_{n-1,j}^{20V}=\sum_{k=1}^{2n-3} 2k \,Z_{n-1,k}^{20V}$. 

The case (c) corresponds 
to the contributions to the refined partition function $Z_{n,2}^{20V \backslash}$ such that a path goes vertically through the left vertex. We must therefore subtract from $Z_{n,2}^{20V \backslash}$ the contributions with the two rightmost paths ending diagonally before their last vertical step, in bijection with 20V configurations on $\cQ_{n-1}$. This gives a net contribution of
$Z_{n,2}^{20V \backslash}-Z_{n-1}^{20V}$. The case (d) contributes
$Z_{n,2}^{20V -}$. Finally the case (e) contributes $Z^{20V}({\mathcal P}_{n,n-3})$. Collecting all the terms and using
$Z_{n,2}^{20V \backslash}+Z_{n,2}^{20V -}=Z_{n,2}^{20V}$, we get:
$$Z^{20V}({\mathcal P}_{n,n-2})=\sum_{k=1}^{2n-3} 2k \,Z_{n-1,k}^{20V}
+Z_{n,2}^{20V}-Z_{n-1}^{20V}+Z^{20V}({\mathcal P}_{n,n-3})$$
Using Remarks \ref{palinrem} and \ref{reciprem}, which imply the symmetry $Z_{n-1,2n-2-k}^{20V}=Z_{n-1,k}^{20V}$, we further note that 
$$\sum_{k=1}^{2n-3} 2k \,Z_{n-1,k}^{20V}=\frac{1}{2}\sum_{k=1}^{2n-3} 2k (Z_{n-1,k}^{20V}+Z_{n-1,2n-2-k}^{20V}) =2(n-1)\sum_{k=1}^{2n-3}Z_{n-1,k}^{20V}=2(n-1)Z_{n-1}^{20V} $$ 
leading to the identity:
\begin{equation}\label{firstrela} Z^{20V}({\mathcal P}_{n,n-3})= Z^{20V}_n-Z^{20V}_{n,2}-2(n-1)Z^{20V}_{n-1}
\end{equation}
For the domino tiling model, we start from $Z^{DT}({\mathcal T}_{n,n-2})$. It splits again according to the configurations of it two topmost vertices, as shown in Figs.\ref{fig:pachter3} (a) and (b). The case (a) corresponds to a path visiting the leftmost vertex, which must be made of $n-1$ first up steps, and continue on the set of dashed edges until it reaches the endpoint $(2n-1,2n-1)$, giving rise to $2(n-1)$ possibilities, while the rest of the configuration is arbitrary on $\cQ_{n-1}$. The corresponding contribution is $2(n-1)Z_{n-1}^{DT}$. In the case (b) a path visits the right vertex, while the left vertex is empty. The corresponding configurations are those paths contributing to the refined partition function $Z_{n,n-2}^{DT}$ that do not pass through the left vertex. We must therefore subtract from $Z_{n,n-2}^{DT}$ the two contributions of the topmost path that visit both vertices in $\cT_{n,n-2}$ (with either a horizontal step between the two vertices, or a down followed by an up step). In both cases, the rest of the paths correspond to an arbitrary tiling of $\cT_{n-1}$. The total contribution is therefore $Z_{n,n-2}^{DT}-2 Z_{n,n-1}^{DT}$. The case (c) with the two empty vertices contributes
$Z^{DT}({\mathcal T}_{n,n-3})$. Assembling all the terms, we get
$$ Z^{DT}({\mathcal T}_{n,n-2})=2(n-1)Z_{n-1}^{DT}+Z_{n,n-2}^{DT}-2 Z_{n,n-1}^{DT}+Z^{DT}({\mathcal T}_{n,n-3}) $$
We deduce that
\begin{equation}\label{secondrela}Z^{DT}({\mathcal T}_{n,n-3})= Z^{DT}_n-Z_{n,n-2}^{DT}-(2n-3)Z^{DT}_{n-1}
\end{equation}
Finally we use eq.\eqref{usefuldt} of Corollary \ref{ref20dt} for $k=n-2$ to identify:  
$Z_{n,n-2}^{DT}=Z_{n,1}^{20V}+Z_{n,2}^{20V}$. Together with the identities $Z^{20V}_m=Z^{DT}_m$ of Theorem \ref{equivalencethm} and the fact that $Z_{n,1}^{20V}=Z_{n-1}^{20V}$ \eqref{obs}, this implies the equality between the r.h.s. of \eqref{firstrela} and \eqref{secondrela}, and the theorem follows.
\end{proof}

%
%
%
%

\section{Discussion and Conclusion}
\label{concsec}

In this paper we have investigated the 20V model with DWBC3 boundary conditions on the family of pentagons $\cP_{n,n-k}$ and their relation to the domino tilings of the family of Aztec triangles $\cT_{n,n-k}$, and proved the conjectured identity between their partition functions for the cases $k=1$, $2$, $3$. In the case $k=1$ of the quadrangle $\cQ_n=\cP_{n,n-1}$ we have extended the result to refined partition functions. Our proofs for this case have relied mainly on the choice of integrable weights for the 20V model, allowing for transforming the original problem into a 6V model with U-turn boundaries. The further truncations of the domain $\cQ_n$ for $k=2,3$ were treated in a purely combinatorial manner, by relating the partition functions to the $k=1$ case and its refinements. 

We now discuss a conjecture for the exact number $Z_n^{20V}=Z_n^{DT}$ and alternative formulas for these quantities. We also show how to introduce some extra weight in the domino tiling problems, before making some concluding remarks.

\subsection{An exact (conjectured) formula}
\label{secconj}

Based on numerical data, we were able to formulate the following conjecture about the total number of configurations of the 20V model on the quadrangle $\cQ_n$:
\begin{conj}\label{numconj}
The total numbers $Z^{20V}_n=Z_n^{DT}$ of configurations of the 20V model on the quadrangle $\cQ_n$ and of domino tilings of the Aztec triangle $\cT_n$ read: 
\begin{equation}\label{conjnum}
Z^{20V}_n=Z_n^{DT}= 2^{n(n-1)/2}\prod_{i=0}^{n-1}  \frac{(4i+2)!}{(n+2i+1)!}
\end{equation}
\end{conj}

This formula was checked up to $n=30$. It is reminiscent of the famous formula for the number $A_n$ of $n\times n$ Alternating Sign Matrices (ASM), which reads
$A_n=\prod_{i=0}^{n-1}\frac{(3i+1)!}{(n+i)!}$
or the related formula for the dimension $\delta_n=3^{n(n-1)/2}\,A_n$ of the irreducible representation of $GL_n$ indexed by the Young diagram $Y_n=(n-1,n-1,n-2,n-2,...,1,1,0,0)$.

Note that the technique used by Kuperberg \cite{kuperberg1996another,kuperberg2002symmetry} for proving the ASM formula and its variants doesn't seem to apply here, as we have not been able to find a factorized formula for the partition function in which the spectral parameters approach their combinatorial value as $z_i=q^4 s^i,w_j=q^{-8}s^{j+n}$ for some parameter $s$ to be sent to $1$ eventually. 

Another approach would try to reproduce the connection between the 6V model with DWBC conditions and the $GL_n$ character with partition $Y_n$ or that between U-turn 6V model and characters of $SP_{2n}$ \cite{stroIK,RS}. Looking for candidates, and denoting by $B_n$ the r.h.s. of \eqref{conjnum}, we have found\footnote{All these results are product formulas obtained by specializing the corresponding characters, given by various determinant formulas \cite{Mizu}, when all arguments are $1$.} that the irreducible representation of $SP_{6n}$ with diagram
$Y_n'=(n-1,n-1,n-1,n-2,n-2,n-2,...,1,1,1,0,0,0)$ has the dimension $\delta_n'=B_{2n}/2^{2n}$. Analogously, we found that 
the irreducible representation of $SO_{6n+1}$ with diagram $Y_n''=(\frac{2n-1}{2},\frac{2n-1}{2},\frac{2n-1}{2},...,\frac{3}{2},\frac{3}{2},\frac{3}{2},\frac{1}{2},\frac{1}{2},\frac{1}{2})$ has dimension $\delta_n''=2^n \, B_{2n}$. Finally, we found that 
the irreducible representation of $SO_{6n-2}$ with diagram $Y_n'''=(n-1,n-1,n-1,n-2,n-2,n-2,...,1,1,1,0,0)$ has dimension 
$\delta_n'''=2^{n-1} \,B_{2n-1}$. Recalling that the in 6V case \cite{stroIK,RS} the argument involved some non-trivial $\Z_3$ symmetry of the partition function at the combinatorial point, we expect a $\Z_4$ symmetry to play an analogous role here. We intend to return to this question in a future publication.

\subsection{Alternative expressions for $Z_n^{20V}$}
\label{altsec}

The results of this paper allow us to derive several alternative determinant formulas for the number $Z_n^{20V}$. 
Using Theorem \ref{thm20V}, we may express $Z_n^{20V}$ as follows.
\begin{thm}\label{binothm}
The number $Z_n^{20V}$ of configurations of the $20V$ model on the quadrangle $\cQ_n$ reads:
\begin{equation}Z_n^{20V} =\det_{0\leq i \leq j\leq n-1} \left( 2^i \, {i+2j+1\choose 2j+1} -{i-1\choose 2j+1} \right) \end{equation}
with the convention ${m\choose p}=0$ for all $-1\leq m<p$.
\end{thm}
\begin{proof}
We start from the determinant formula $Z_n^{20V}=\det_{0\leq i,j \leq n-1}(f_P(u,v)\vert_{u^iv^j})=\det(P_n)$ of Theorem \ref{thm20V}, where:
$$f_P(u,v)=g_{20V}(u,v)=
\frac{\left(\frac{1+u}{1-u}\right)^{2}}{1-v\,\left(\frac{1+u}{1-u}\right)^{2}} -\frac{u^{2}}{1-v\, u^{2}}
$$
and $P_n$ is the finite truncation of the infinite matrix $P$ generated by $f_P$.
Using the same technique as in Sections \ref{detrelasec}-\ref{equivsec}, let us consider the infinite matrix $\Lambda$ with generating function:
$$f_\Lambda(u,v)=\frac{1-u}{1-u-vu} $$
Then (1) $\Lambda$ is lower triangular and (2) the diagonal matrix elements of $\Lambda$ are all $1$. We deduce 
the truncation $(\Lambda\, P)_n=\Lambda_n\,P_n$, and therefore $\det((\Lambda\, P)_n)=\det(P_n)$. Explicitly computing:
$$f_{\Lambda\, P}(u,v)=(f_\Lambda*f_P)(u,v)=f_P\left(\frac{u}{1-u},v\right) =
\frac{\left(\frac{1}{1-2u}\right)^{2}}{1-v\,\left(\frac{1}{1-2u}\right)^{2}} -\frac{\left(\frac{u}{1-u}\right)^{2}}{1-v\, \left(\frac{u}{1-u}\right)^2}$$
we identify the entries
$$(\Lambda P)_{i,j}=\left\{\left(\frac{1}{1-2u}\right)^{2j+2}-\left(\frac{u}{1-u}\right)^{2j+2}\right\}\Bigg\vert_{u^i}=2^i \, {i+2j+1\choose 2j+1} -{i-1\choose 2j+1}$$
and the theorem follows.
\end{proof}

A slight variant of the formula consists of using the analytic continuation of binomial coefficients.
\begin{cor}\label{symcor}
The number $Z_n^{20V}$ of configurations of the $20V$ model on the quadrangle $\cQ_n$ reads:
\begin{equation}\label{slight}
Z_n^{20V}=\frac{2^{n(n-1)/4}}{2}\, \det_{0\leq i,j\leq n-1}\Big(\theta_{2j+1}(i)+\theta_{2j+1}(-i)  \Big), \quad \theta_m(x)=\frac{2^{x/2}}{m!} (x+1)(x+2)...(x+m)
\end{equation}
\end{cor}

This latter, more symmetric form allows to derive the following constant term identity.
\begin{thm}\label{ctthm}
The number $Z_n^{20V}$ of configurations of the $20V$ model on the quadrangle $\cQ_n$ reads:
$$Z_n^{20V}=CT_{x_1,...,x_n}\left\{ 
\frac{\prod_{1\leq i<j\leq n} (x_i-x_j)(1+x_i+x_j-x_ix_j)}{\prod_{i=1}^n x_i^{2i-1}\, (1-x_i)^n}  \right\} $$
where the symbol $CT_{x_1,...,x_n}$ stands for the constant coefficient in the Laurent series expansion around $x_i=0$ for all $i$.
\end{thm}
\begin{proof}
Let us rewrite the function $\theta_m(x)$ of \eqref{slight} as:
$$\theta_m(x) =\frac{1}{m!}\partial_u^m \, u^{m+x}\big\vert_{u=\sqrt{2}} =CT_t \left\{ e^{t\,\partial_u}  t^{-m} \, u^{m+x}\right\}\big\vert_{u=\sqrt{2}}=CT_t \left\{   t^{-m} \, (\sqrt{2}+t)^{m+x}\right\}$$
where we have used the Taylor expansion to interpret the operator $e^{t\,\partial_u}$ as
the shift map $S_{t}: f(u)\mapsto f(u+t)$. We now use multilinearity of determinants to rewrite
$$
Z_n^{20V}=\frac{2^{n(n-1)/4}}{2}\, CT_{t_1,t_2,...,t_n}\left\{
\prod_{j=1}^n \left(1+\frac{\sqrt{2}}{t_j}\right)^{2j-1} \, \det_{1\leq i,j\leq n}\left( (\sqrt{2}+t_j)^{i-1}+(\sqrt{2}+t_j)^{1-i}\right)\right\}
$$
Next we use the standard formula
$$ \det_{1\leq i,j\leq n}( x_j^{i-1}+x_j^{1-i})= 2\,\frac{\prod_{1\leq i<j\leq n} (x_i-x_j)(1-x_ix_j)}{\prod_{i=1}^n x_i^{n-1}}$$
to rewrite:
\begin{equation}\label{ctsym}
Z_n^{20V}=2^{n(n-1)/4}\,CT_{t_1,t_2,...,t_n}\left\{
\prod_{i=1}^n \left(1+\frac{\sqrt{2}}{t_i}\right)^{2i-1} \, \prod_{1\leq i<j\leq n} (t_j-t_i)\left(1-\frac{1}{(\sqrt{2}+t_i)(\sqrt{2}+t_j)}\right)
\right\}
\end{equation}
We now use the obvious property of the constant term that $CT_\bx(f(\bx))=CT_\bx({\rm Sym}(f(\bx))$, where for any function 
of the variables $\bx=(x_1,...,x_n)$ we define the symmetrization map 
$${\rm Sym}: f(\bx)\mapsto {\rm Sym}(f)(\bx):=\frac{1}{n!} \sum_{\sigma\in S_n} f(x_{\sigma(1)},...,x_{\sigma(n)}) .
$$
We also need the antisymmetrization map
$${\rm ASym}: f(\bx)\mapsto {\rm ASym}(f)(\bx):=\frac{1}{n!} \sum_{\sigma\in S_n} {\rm sgn}(\sigma)\,f(x_{\sigma(1)},...,x_{\sigma(n)})
$$
with the main property
$$
{\rm Sym}\left( \prod_{1\leq i<j\leq n}(x_j-x_i) f(x_1,...,x_n) \right)=\prod_{1\leq i<j\leq n}(x_j-x_i)\,{\rm ASym}\left(f(x_1,...,x_n)\right) .
$$
Symmetrizing the argument of the constant term in  \eqref{ctsym}, we get:
\begin{eqnarray*}
Z_n^{20V}&=&\frac{2^{n(n-1)/2}}{n!}CT_{t_1,t_2,...,t_n}\left\{ \prod_{i=1}^n \left(1+\frac{\sqrt{2}}{t_i}\right) \right. \\
&&\quad \times \, \left. \prod_{1\leq i<j\leq n}\left(1-\frac{1}{(\sqrt{2}+t_i)(\sqrt{2}+t_j)}\right)(t_j-t_i)\left(\left(1+\frac{\sqrt{2}}{t_j}\right)^2-\left(1+\frac{\sqrt{2}}{t_i}\right)^2\right)\right\} .
\end{eqnarray*}
Let us now change variables to $x_i=\frac{t_i}{t_i+\sqrt{2}}$. Viewing the constant term in $t_j$ as a contour integral
$\oint \frac{dt_j}{2\pi i \,t_j}=\oint \frac{dx_j}{2\pi i \, x_j(1-x_j)}$, we easily find:
\begin{eqnarray*}
Z_n^{20V}&=&\frac{1}{n!}CT_{x_1,x_2,...,x_n}\left\{
\frac{\prod_{1\leq i<j\leq n}(1+x_i+x_j-x_ix_j)(x_i-x_j)(x_j^{-2}-x_i^{-2})}{ \prod_{i=1}^n x_i\,(1-x_i)^n}\right\}\\
&=& CT_{x_1,x_2,...,x_n}\left\{{\rm Sym}\left(
\frac{\prod_{1\leq i<j\leq n}(1+x_i+x_j-x_ix_j)(x_i-x_j)}{ \prod_{i=1}^n x_i^{2i-1} \,(1-x_i)^n}\right)\right\} .
\end{eqnarray*}
The theorem follows.
\end{proof}

The constant term identity of Theorem \ref{ctthm} is reminiscent of constant term expressions for TSSCPP derived by Zeilberger \cite{zeilberger1}  in his famous proof of the ASM conjecture.

\subsection{Introducing step weights in the domino tiling problem}
\label{stepweightsec}

In this paper we have found compact determinant formulas for both the number of tilings of the Aztec triangle $\cT_n$ and its refinements, involving the truncation of infinite matrices generated  respectively by the series $f_{DT}$ \eqref{fat} and $f_{DT}^{\rm ref}$ \eqref{refat}.
We note that a very simple decoration of the generating function $f_{DT}(u,v)$ allows to include a non-trivial extra step weight for the non-intersecting Schr\"oder paths. In general, one would want to associate weights $\al,\beta,\gamma$ respectively to up, down and horizontal steps. However, by simple rescaling of $u,v$ we may restrict without loss of generality to only a weight $\gamma$ per horizontal step.

\begin{thm}\label{thmATgamma}
The partition function for domino tilings of the Aztec triangle of order $n$, with an extra weight $\gamma$ per horizontal step in the non-intersecting Schr\"oder path formulation reads:
\begin{eqnarray}
Z^{DT,\gamma}_n&=&\det_{0\leq i,j \leq n-1}\left( 
f_{DT,\gamma}(u,v)\Big\vert_{u^iv^j} \right)\label{zatgamma} \\
f_{DT,\gamma}(u,v)&=& \frac{1+u}{1-v -2(1+\gamma) u v - u^2 v +\gamma^2 u^2 v^2} \label{fatgamma}
\end{eqnarray}
\end{thm}
\begin{proof}
We adapt the proof of Theorem \ref{genDT}. By the Lindstr\"om-Gessel-Viennot theorem, the partition function function of the non-intersecting Schr\"oder paths with the new weights is given by $Z^{DT,\gamma}_n=\det((M_\gamma)_n)$, where the subscript $n$ indicates the truncation to the first $n$ rows and columns of the infinite matrix $M_\gamma$
whose entries $(M_\gamma)_{i,j}$ are the partition functions for a single path form $(0,2i)$ to $(2j+1,2j+1)$,
counted with a multiplicative weight $\gamma$ per horizontal step. We simply have to show that these entries are generated by:
$$f_{M_\gamma}(u,v)=\frac{1+u}{1-v -2(1+\gamma) u v - u^2 v +\gamma^2 u^2 v^2}=f_{DT,\gamma}(u,v)$$
To show this, we start with the function $\sigma_\gamma(u,v):=\frac{1}{1-u-v-\gamma u v}$ the partition function for arbitrary finite length Schr\"oder paths from the origin, with weight $u$ per up step, $v$ per down step, and $\gamma uv$ per horizontal step. Accordingly the partition function $(\Sigma_\gamma)_{i,j}$ of paths from the origin
with fixed endpoint $(i,j)$, with $i-j$ even, is given by
$(\Sigma_\gamma)_{i,j}=\sigma_\gamma(u,v)\vert_{u^{\frac{i+j}{2}} v^{\frac{i-j}{2}}}$, or equivalently by a straightforward generalization of \eqref{sigentry}. Noting again that $(M_\gamma)_{i,j}=(\Sigma_\gamma)_{2j+1,2j+1-2i}$, we get
$$(M_\gamma)_{i,j}=\oint \frac{dx}{2i\pi x} \frac{dy}{2i\pi y}\, \frac{1}{x^{2j+1-i} y^{i}} \, \frac{1}{1-x-y-\gamma x y}=\oint \frac{du}{2i\pi u} \frac{dv}{2i\pi v}\, \frac{1}{v^{2j+1} u^{i}} \, \frac{1}{1-v-u v-\gamma u v^2}$$
after again changing integration variables to $u=y/x, v=x$.
Using the same trick \eqref{trick} as in previous section, we conclude that:
$$f_{M_\gamma}(u,v)=\frac{\frac{1}{1-(1+u)\sqrt{v}-\gamma u v}-\frac{1}{1+(1+u)\sqrt{v}-\gamma u v}}{2\sqrt{v}}=\frac{1+u}{(1-\gamma u v)^2-v(1+u)^2}$$
and the Theorem follows.
\end{proof}

The result above is easily extended to the refined domino tiling of Section. \ref{refdtsec}.
The partition function $Z_n^{DT,\gamma}(t):=\sum_{k=1}^n Z^{DT,\gamma}_{n,k} \,t^{k-1}$ now including a multiplicative weight $\gamma$ per horizontal step in the non-intersecting Schr\"oder path formulation of the refined tiling problem, can be directly calculated as follows (we leave the proof to the reader as a straightforward exercise).

\begin{thm}\label{refrefdominothm}
The generating function $Z^{DT}(t,\gamma)$ reads:
\begin{eqnarray*}Z^{DT}(t,\gamma)&=&\det_{0\leq i,j\leq n-1} \left(
f_{DT,\gamma}^{\rm ref}(u,v)\Big\vert_{u^iv^j}\right) \\
f_{DT,\gamma}^{\rm ref}(u,v)&=&f_{DT,\gamma}(u,v)+v^{n-1} \frac{\al_+(\gamma;u)^{2n}}{\sqrt{1+2(1+2\gamma)u+u^2}}\, \frac{(t-1)u}{\al_+(\gamma;u)-t u}
\end{eqnarray*}
where $\al_{+}(\gamma,u)= \frac{1}{2}\big(1 + u+\sqrt{1 + 2(1+2\gamma) u + u^2}\big)$ and $f_{DT,\gamma}(u,v)$ as in \eqref{fatgamma}.
\end{thm}

It is interesting to notice that the upper triangular matrix $L$ involved in the proofs of Theorems \ref{equivalencethm} and \ref{fullequivthm} can be modified to obtain an alternative expression for the partition function generating function  $Z^{DT,\gamma}(t)$ 
above. 
Indeed, let $L_\gamma$ be the infinite matrix generated by:
$$f_{L_\gamma}(u,v):= \frac{1+2\gamma u-\gamma u^2}{1-u-v u(1+\gamma u)}$$
Writing $Z^{DT}(t,\gamma)=\det\big( (L_\gamma {\tilde M}_\gamma^{(n)}(t))_n\big)$, where $ {\tilde M}_\gamma^{(n)}(t)$ is generated by $f_{DT,\gamma}^{\rm ref}(u,v)$ above, we finally get:
\begin{thm}
The generating function $Z^{DT}(t,\gamma)$ reads:
\begin{eqnarray*}
Z^{DT}(t,\gamma)&=& \det\Big( 
{\bar f}_{DT,\gamma}^{\rm ref}(u,v)\big\vert_{u^iv^j}
\Big)\\
{\bar f}_{DT,\gamma}^{\rm ref}(u,v)
&=&\frac{(1+\gamma u^2)(1+2\gamma u-\gamma u^2)}{(1-\gamma^2u^2v)((1-u)^2-v(1+\gamma u)^2)}+v^{n-1} \left(\frac{1+\gamma u}{1-u}\right)^{2n} \frac{(t-1)u}{1- t u} 
\end{eqnarray*}
\end{thm}
\begin{proof}
We use the infinite matrix ${\bar M}_\gamma^{(n)}(t):=L_\gamma {\tilde M}_\gamma^{(n)}(t)$, whose generating function is easily computed via the formula 
$f_{L_\gamma \,{\tilde M}_\gamma^{(n)}(t)}(u,v)
=\frac{1+2\gamma u-\gamma u^2}{1-u}f_{DT,\gamma}^{\rm ref}\left(u\frac{1+\gamma u}{1-u},v\right)={\bar f}_{DT,\gamma}^{\rm ref}(u,v)$.
\end{proof}
Note also that 
$$\frac{(1+\gamma u^2)(1+2\gamma u-\gamma u^2)}{(1-\gamma^2u^2v)((1-u)^2-v(1+\gamma u)^2)}=\frac{\left(\frac{1+\gamma u}{1-u}\right)^2}{1-v\left(\frac{1+\gamma u}{1-u}\right)^2 } -\frac{\gamma^2 u^2}{1-v \gamma^2 u^2}$$
This $\gamma$-deformation of \eqref{f20v}  is very suggestive, however we have not been able to find a counterpart of the parameter $\gamma$ in the 20V model.

\subsection{Conclusion}
\label{seconc}

Despite the great progress on proving the conjectures of Ref. \cite{DG19}, we are still only understanding part of the sequence of correspondences between the 20V model on the domains $\cP_{n,k}$ and the domino tiling of the domains $\cT_{n,k}$. One way to gain a better global understanding of these would be to include deformation (spectral?) parameters in the correspondence of partition functions on $\cQ_n$ and $\cT_n$ respectively. However, the number of deformation parameters at hand is limited. For instance, the trick we used with the last $w_n=w$ vertical spectral parameter cannot be adapted to horizontal spectral parameters: these are indeed all frozen to the same value $z=p^{-2}q^{-2}$ in order for the relation to the U-turn boundary 6V partition function to hold. One may however try to keep more non-trivial values of the vertical spectral parameters $w_n,w_{n-1},...$ etc.

Another direction of study regards large $n$ asymptotics. We expect the typical large configurations of the 20V model on $\cQ_n$ to exhibit an arctic curve phenomenon, such as that observed in \cite{DDG20} for the 20V DWBC1,2 models, and derived analytically using the tangent method of Colomo and Sportiello \cite{COSPO}.  We intend to address this question in a future publication.

\bibliographystyle{amsalpha}
\bibliography{20VDW3}
\end{document}